\documentclass{IEEEoj}
\usepackage{cite}
\usepackage{amsmath,amssymb,amsfonts}
\usepackage{algorithmic}
\usepackage{graphicx,color}
\usepackage{textcomp}
\usepackage{listings}
\usepackage{xfrac}
\usepackage{url}
\def\BibTeX{{\rm B\kern-.05em{\sc i\kern-.025em b}\kern-.08em
    T\kern-.1667em\lower.7ex\hbox{E}\kern-.125emX}}
\AtBeginDocument{\definecolor{ojcolor}{rgb}{.098039,0.003922,0.71765}}
\def\OJlogo{\includegraphics[height=20pt]{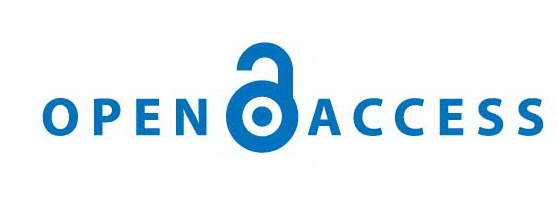}}
\def\OJlogoii{}
\makeatletter
\def\ps@headings{%
  \def\@oddhead{\vbox{\hbox to \textwidth{\OJlogoii\hfill\OJlogo}\par
  \vspace*{0pt}\hbox to \textwidth{\vrule width\textwidth height.3pt depth0pt}}}%
 \def\@evenhead{\vbox{\hsize\textwidth\vbox to 0pt{\hsize\textwidth\vspace*{7.7pt}\rfxfont\raggedright\rightmark:\ \leftmark\hfill\par}\par\vspace*{16pt}\hbox to \textwidth{\vrule width\textwidth height.3pt depth0pt}}}%
        \def\@evenfoot{\hbox to \textwidth{{\rffont\thepage}\hfill{{\rffont VOLUME\ \@jvol,\ \@pubyear}}}}%
        \def\@oddfoot{\hbox to \textwidth{{{\rffont VOLUME\ \@jvol,\ \@pubyear}}\hfill{\rffont\thepage}}}%
	          }%
\def\ps@plain{%
  \def\@oddhead{\vbox{\hbox to \textwidth{\OJlogo\hfill\OJlogoii}\par
  \vspace*{0pt}\hbox to \textwidth{\vrule width\textwidth height.3pt depth0pt}}}%
   \let\@evenhead\@oddhead%
        \def\@evenfoot{\vbox to 10pt{\hbox to \textwidth{\hfill\rffont This work is licensed under a Creative Commons Attribution 4.0 License. For more information, see https://creativecommons.org/licenses/by/4.0/\hfill}\par\vspace*{-12pt}%
                        \hbox to \textwidth{{\rffont\thepage}\hfill{{\rffont VOLUME\ \@jvol,\ \@pubyear}}}}}%
        \def\@oddfoot{\vbox to 10pt{\hbox to \textwidth{\hfill\rffont This work is licensed under a Creative Commons Attribution 4.0 License. For more information, see https://creativecommons.org/licenses/by/4.0/\hfill}\par\vspace*{-12pt}%
                        \hbox to \textwidth{{{\rffont VOLUME\ \@jvol,\ \@pubyear}}\hfill{\rffont\thepage}}}}%
	          }%

\makeatother

\begin{document}
\receiveddate{30 January, 2024}
\reviseddate{30 March, 2024}
%\accepteddate{30 March, 2024}
\publisheddate{30 March, 2024}
%\currentdate{30 March, 2024}
%\doiinfo{OAJPE.2020.2976889}

\title{Chebyshev and The Fast Fourier Transform Methods for Signal Interpolation}

\author{Ishmael N. Amartey\authorrefmark{1}}
\affil{Department of Statistics and Actuarial Science,  Northern Illinois University, Dekalb, IL, 60115, USA}
%\affil{Department of Physics, Colorado State University, Fort Collins, 
%CO 80523 USA}
\corresp{CORRESPONDING AUTHOR: Ishmael N. Amartey (e-mail: iamartey1@niu.edu).}
%\authornote{This work was supported by the Natural Sciences and Engineering Research Council (NSERC) of Canada.}
%\markboth{Preparation of Manuscripts for IEEE TRANSACTIONS and JOURNALS}{Author \textit{et al.}}
%%%%%%%%%%%%%%%%%%%%%%%%%%%%%%%%%%%%%%%%%%%%%%%%%%%%%%%%%%%%%%%%%%%%%%%%%%%

\begin{abstract}
Approximation theorem is one of the most important aspects of numerical analysis that has evolved over the years with many different approaches. Some of the most popular approximation methods include the Lebesgue approximation theorem, the Weierstrass approximation, and the Fourier approximation theorem. The limitations associated with various approximation methods are too crucial to ignore, and thus, the nature of a specific dataset may require using a specific approximation method for such estimates.
In this report, we shall delve into Chebyshev’s polynomials interpolation in detail as an alternative approach to reconstructing signals and compare the reconstruction to that of the Fourier polynomials. We will also explore the advantages and limitations of the Chebyshev polynomials and discuss in detail their mathematical formulation and equivalence to the cosine function over a given interval [a, b].
\end{abstract}


\maketitle
\section{Introduction}
\label{sec:introduction}
The Chebyshev interpolation and points have been shown to have an advantage over other approximation methods due to their unique ability to absorb two parameters for interpolations, unlike other methods that use only one.  In this study, we investigate the Chebyshev interpolation and its properties under different data scenarios using the Chebfun package in Matlab by Trefethen \cite{doi:10.1137/1.9781611975949}.

In section two of this report, we interpolate the Chebyshev points through random data and compute the rounding errors as well as the computation of the Chebyshev points of the first kind. A geometric mean distance between the points, a convergence of the interpolants, and the scaling of the Chebyshev function to the interval [a, b] is also discussed utilizing exercises from Chapter 2 of Trefethen \cite{doi:10.1137/1.9781611975949}. 

Moving on to the third section, we delve into the Chebyshev polynomials and series. This includes an exploration of the dependency on wave numbers, the representation of complex functions using the Chebyshev series, the conditioning of the Chebyshev basis, and an examination of the extrema and roots of Chebyshev polynomials.

Finally, in section 4 an interpolation of the gamma variate function is undertaken utilizing Chebyshev polynomials. Two distinct approaches are employed: 1) employing unevenly distributed Chebyshev nodes, and 2) utilizing evenly distributed Chebyshev nodes. The outcomes derived from these two approaches are compared with those obtained through Fourier polynomials, leading to the formulation of pertinent conclusions, and facilitating in-depth discussions.

%%%%%%%%%%%%%%%%%%%%%%%%%%%%%%%%%%%%%%%%%%%%%%%%%%%

\section{Chebyshev points and interpolants }
The Chebyshev points are scaled on the interval $[-1, 1]$, so for $n+1$ equally spaced angles $\{\theta_j\}$ from $0$ to $\pi$, the Chebyshev points can be thought of as the real parts of $\{z_j\}$ for $n+1$ points on the upper half of the unit circle in the complex plane for $z_j = e^{-it}$. Thus,
\begin{align}
x_j &= \operatorname{Re}(z_j) = \frac{1}{2}(z_j + z_j^{-1}), & 0 &\leq j \leq n \notag \\
\end{align}

The Chebyshev points in their original angle are defined as
\begin{align}
x_j &= \cos{\left(\frac{j\pi}{n}\right)}, & 0 &\leq j \leq n \notag
\end{align}

But this lacks symmetry. To guarantee symmetry, the alternate definition of the Chebyshev points is set to
\begin{align}
x_j &= \cos\left(\frac{(2j+1)\pi}{2n}\right) & \text{for } j &= 0, 1, \ldots, n-1 \notag
\end{align}

Since the term $(2j + 1)$ inside the cosine function is always an odd integer for any integer $j$, the angles inside the cosine function are always odd multiples of $\pi/2n$, and since the cosine function is an even function ($\cos(-\theta) = \cos(\theta)$), this ensures that the Chebyshev points are perfectly symmetric around the origin (i.e., symmetric with respect to $x = 0$) for all values of $n$.

In floating-point arithmetic, this symmetry is preserved because the cosine function is well-behaved for small angles. Mathematically, the equivalence of the Chebyshev points as a cosine function is established below.
\begin{align}
x_j &= \cos{\left(\frac{j\pi}{n}\right)}
\end{align}

Set
\begin{align}
\frac{(2j+1)\pi}{2n} &= \frac{(2\pi+\pi)}{2n} = \frac{(j\pi+\pi/2)}{n}
\end{align}

Rewriting Eq 6 gives
\begin{align}
x_j &= \cos\left(\frac{(j\pi+\pi/2)}{n}\right)
\end{align}

By the trigonometric identity for the cosine function,
\begin{align}
\cos{\left(\theta + \frac{\pi}{2}\right)} &= -\sin{(\theta)}
\end{align}

So
\begin{align}
x_j &= -\sin{(j\pi/n)}
\end{align}

Since
\begin{align}
\sin{\left(\theta\right)} &= \cos{(\pi/2-\theta)}
\end{align}

Eq 10
\begin{align}
x_j &= \cos{\left(\frac{\pi}{2}-\frac{j\pi}{n}\right)}
\end{align}

So
\begin{align}
\cos{\left(\frac{\pi}{2}-\frac{j\pi}{n}\right)} &\approx \cos{\left(\frac{j\pi}{n}\right)}
\end{align}

Therefore
\begin{align}
x_j &= \frac{(2j+1)\pi}{2n} = \cos{\left(\frac{j\pi}{n}\right)}
\end{align}

\begin{align}
    x_j=\cos{\left(\frac{\pi}{2}-\frac{j\pi}{n}\right)}
\end{align}
So,
\begin{align}
    \cos{\left(\frac{\pi}{2}-\frac{j\pi}{n}\right)}\ \approx \cos{\left(\frac{j\pi}{n}\right)}
\end{align}
Therefore
\begin{align}
    x_j=\frac{(2j+1)\pi}{2n}\ = \cos{\left(\frac{j\pi}{n}\right)}
\end{align}

\subsection{Chebyshev Interpolation through Random Data}
Chebyshev polynomials are excellent in interpolating equally spaced points compared to other polynomial interpolations which does not do a great job. This is mainly because of the cluster of points along the ends of the interval. Additionally, the Chebyshev interpolant is effective because it has the same average distance from each point.
In Figures \ref{fig1}-\ref{fig4} below, the chebfun interpolation plots through random data for 10, 100, 1000, and 10000 points are shown respectively. In each case, the minimum and maximum points are shown using the “minandmax(p)” function along with the computer time required for the computation. The red points along the edges of Figure 2-4 are the points within the interval [0.9999,1].
These plots demonstrate the robustness of the Chebyshev interpolants as the number of points we plot does not have much mathematical difference. As shown in the plots, increasing the number of points produces a messy plot, however, the Chebyshev points are all clustered around -1 and 1.

\begin{figure}[!t]
    \centering
    \includegraphics[width =8cm, height=10cm]{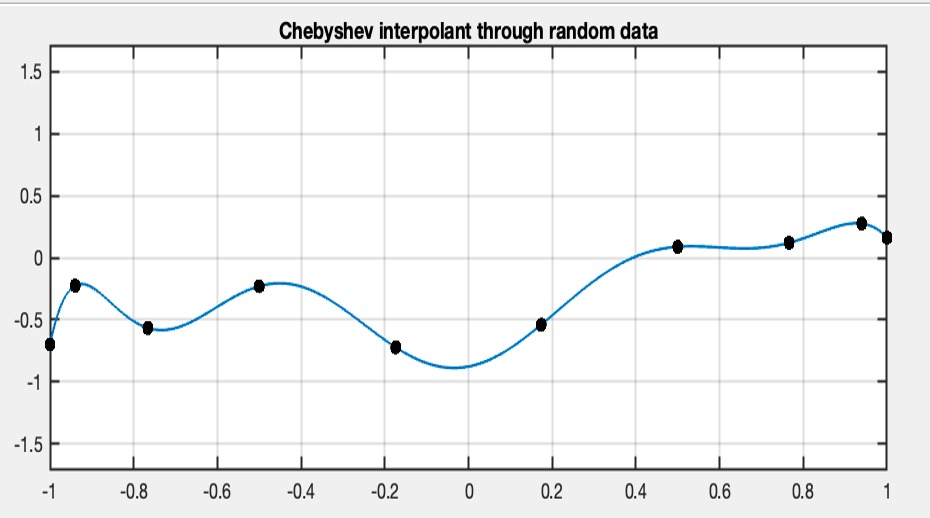}
    \caption{A Chebfun interpolation through random data for 10 points. Elapsed time = 0.28 seconds, min=-0.8299, max=1.1792}
    \label{fig1}
\end{figure}

\begin{figure}[!t]
    \centering
    \includegraphics[width = 8cm, height=10cm]{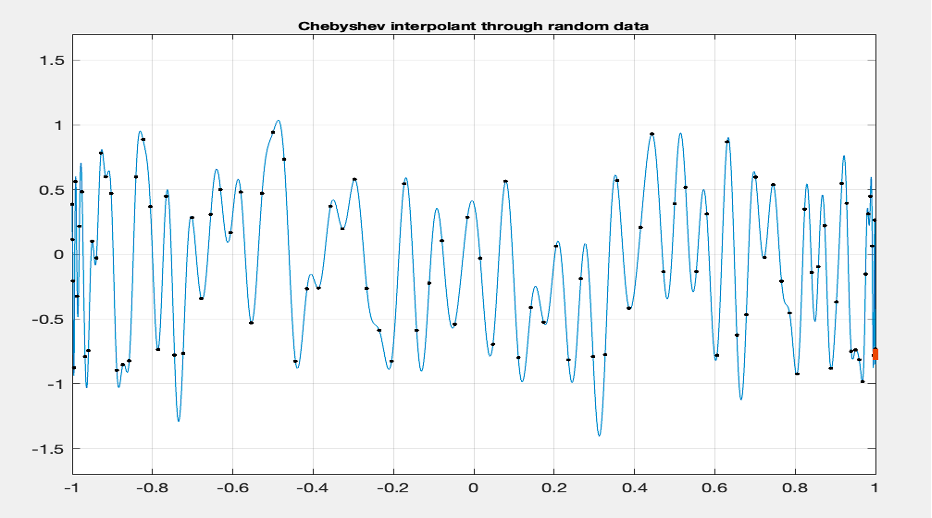}
    \caption{A Chebfun interpolation through random data for 100 points. Elapsed time = 0.07 seconds, min=-1.4093, max = 1.3018.}
    \label{fig2}
\end{figure}

\begin{figure}[!t]
    \centering
    \includegraphics[width = 8cm, height=10cm]{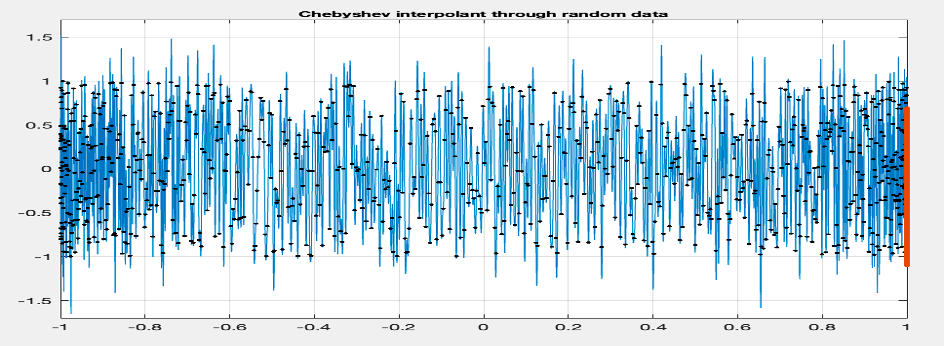}
    \caption{A Chebfun interpolation through random data for 1000 points. Elapsed time = 0.084655 seconds, min = -1.5144, max = 1.8876}
    \label{fig3}
\end{figure}

\begin{figure}[!t]
    \centering
    \includegraphics[width = 8cm, height=10cm]{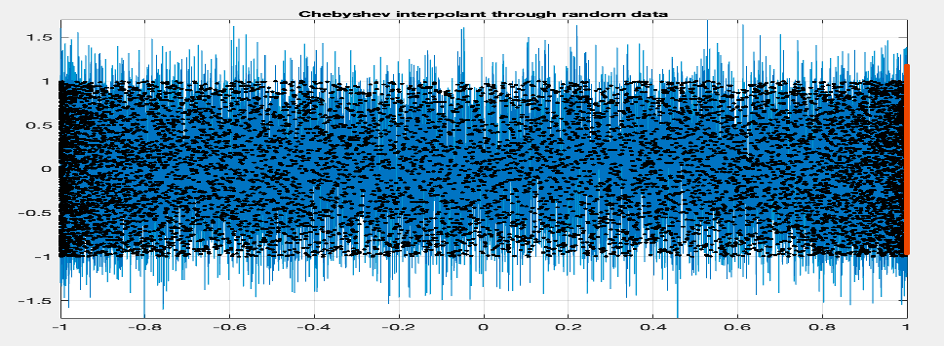}
    \caption{A Chebfun interpolation through random data for 10000 points. Elapsed time = 0.116381 seconds, min =-1.8462, max= 1.7696}
    \label{fig4}
\end{figure}

\subsection{Rounding errors in computing Chebyshev points}

The Matlab program that finds the smallest even value $n \geq 2$ as computed by $x_{\sfrac{n}{2}} \neq 0$ is given below:
\begin{lstlisting}
n = 2; % Start with n = 2
while true
    x = cos((0:n)*pi/n);
    if x(n/2 + 1) ~= 0
        break; 
    end
    n = n + 2; 
end
\end{lstlisting}
disp(['Smallest even n for which x', num2str(n/2), ' $ \neq $ 0: ', num2str(n)]);

\subsection{Chebyshev points of the first kind}

The Chebyshev points, known as Chebyshev points of the first kind or Gauss–Chebyshev points are obtained by taking the real parts of points on the unit circle as already mentioned.  Another function that exhibits similar characteristics of the Chebyshev interpolants is the Legendre points. Just as the Chebyshev points, the Legendre points also have approximately the same average distance between points following the geometric mean estimates.

For the specific case of $n + 1 = 100$, we carry out an investigation into the maximum difference between the Chebyshev point of the first kind and the corresponding Legendre point using the “chebpts” function in the chebfun package in Matlab for the Chebyshev points and the “legpts” function for the Legendre points.

\begin{figure}[!t]
    \centering
    \includegraphics[width = 8cm, height=10cm]{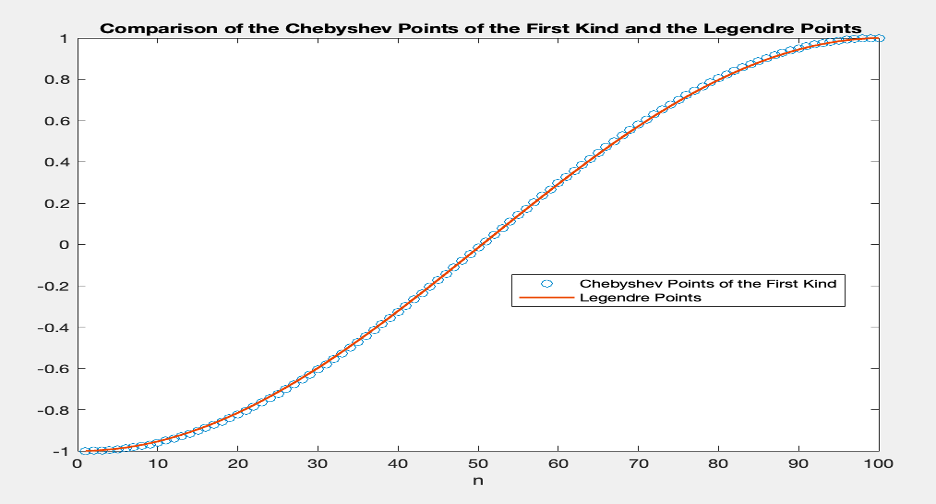}
    \caption{A comparison of the Chebyshev points of the first kind and Legendre points. The maximum between the two points is 0.0084. }
    \label{fig5}
\end{figure}

Figure \ref{fig5}  depicts the illustration of the Chebyshev points and the Legendre points for $n=99$. Clearly, it becomes evident that the Chebyshev points of the first kind and the corresponding Legendre points exhibit a high degree of similarity. The visual representation of the plots reveals a close alignment between the two sets of points with the maximum difference between them being 0.0084.

\subsection{Convergence of Chebyshev interpolants}

Consider a function \( f\left(x\right)=e^x \) on \([-1,1]\) and \( g\left(x\right)=\frac{1}{(1+25x^2)} \). We shall test the convergence of these functions using the ``chebfun'' command on them. But before that, we checked how large \( n \) must be for the level of machine precision using the following commands.

\begin{lstlisting}
machine_epsilon = eps; % Get the machine epsilon value
n_max = 1000;
n = 1;
while n <= n_max
    e1 = norm(f1 - chebfun(fx1, n));
    e2 = norm(f2 - chebfun(fx2, n)); 
    if e1 < machine_epsilon && e2 < machine_epsilon
        fprintf('For machine precision, n must be at least %d\n', n);
        break;
    end
    n = n + 1;
end
if n > n_max
    fprintf('Machine precision not reached within n_max iterations.\n');
end
\end{lstlisting}

For machine precision, \( n \) must be at least 216 so if \( n \) is increased beyond this point, there wouldn’t be any change in accuracy. Using ``chebfun'' the log scale plots of \( \parallel f-p_n\parallel \) were produced with \( p_n \) being the Chebyshev interpolant and \( f \) and \( g \) the functions to be interpolated. We take \( \parallel\bullet\parallel \) to be the supremum norm computed as norm(f-p). Figure 6 depicts the semilog plots for the \( f \) and \( g \). The downward slope followed by a nearly straight line in the log scale for the \( f\left(x\right)=e^x \) suggests exponential convergence. On the other hand, the straight line in the log scale for \( g\left(x\right)=\frac{1}{(1+25x^2)} \) suggests algebraic convergence. These results imply that \( f\left(x\right)=e^x \) minimizes the interpolation error more rapidly as \( n \) increases compared to \( g\left(x\right)=\frac{1}{(1+25x^2)} \).

\begin{figure}[!t]
    \centering
    \includegraphics[width = 8cm, height=10cm]{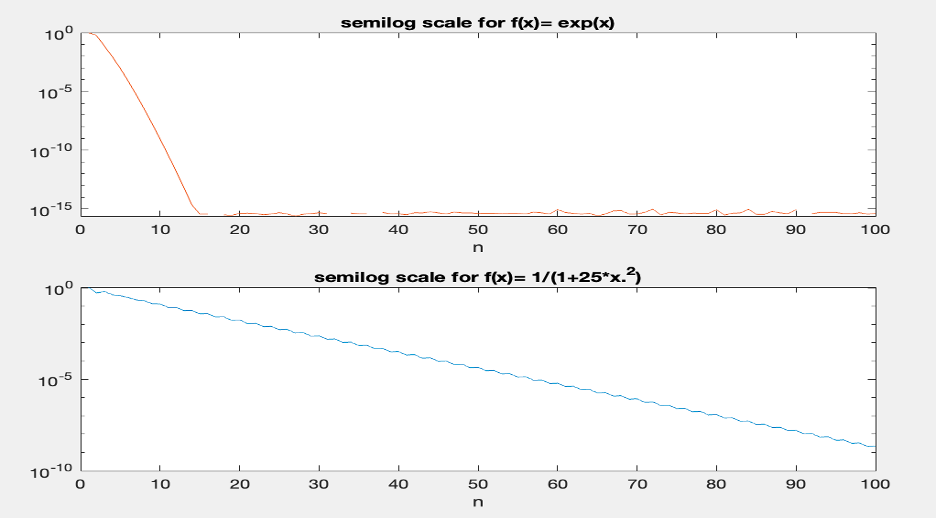}
    \caption{A log scale plot of $\left| f - p_n \right| \text{ as a function of } n \text{ for } f\left(x\right) = e^x \text{ on } [-1,1]$ where $g\left(x\right) = \frac{1}{(1+25x^2)}$}
    \label{fig6}
\end{figure}

\subsection{Geometric mean distance between points}

As mentioned in 2.3, the Chebyshev points have the same average distance between points following the geometric mean estimates. The ‘meandistance’ function created in in Matlab with excepts in the appendix takes a vector of points $x_0......x_n$ within the interval [-1, 1] as input. The code then generates a plot where $x_j$ is plotted on the horizontal axis, and the geometric mean of the distances from $x_j$ to the other points is plotted on the vertical axis utilizing the ‘prod’. function in Matlab. The results are investigated for three different sets of points: Chebyshev points, Legendre points, and equally spaced points within the interval [-1, 1].

\subsubsection{Chebyshev Points}
For $n = 5, 10, 20$, the ‘meandistance’ function code will compute and illustrate the geometric mean distance between each Chebyshev point and the other points in the set providing insights into how the distribution and density of Chebyshev points change as the value of n increases.

\subsubsection{Legendre Points}
Similarly, the code will be employed to analyze the geometric mean distance for Legendre points, as introduced in 2.3. This comparison will highlight any notable similarities or differences between Chebyshev and Legendre point distributions.

\subsubsection{Equally Spaced Points}
Finally, the code will be applied to equally spaced points within the interval [-1, 1]. This comparison serves as a benchmark, allowing us to observe how the geometric mean distance behaves for a straightforward, regularly spaced set of points.

The outcomes of these analyses as shown in Figure \ref{7}- \ref{fig9}, highlight the similarities between the Chebyshev points and that of the Legendre points. For a small n value of 5, the Chebyshev point stretches to its boundary points whereas the Legendre points did not. However, as the number of n increases, both points turn to clustered around their boundary points, but the Chebyshev points are much faster compared to the Legendre points. For the equally spaced points, there isn’t much of a difference as the number of points is spread equally across the length of $x_j$.

\begin{figure}[!t]
    \centering
    \includegraphics[width = 8cm, height=10cm]{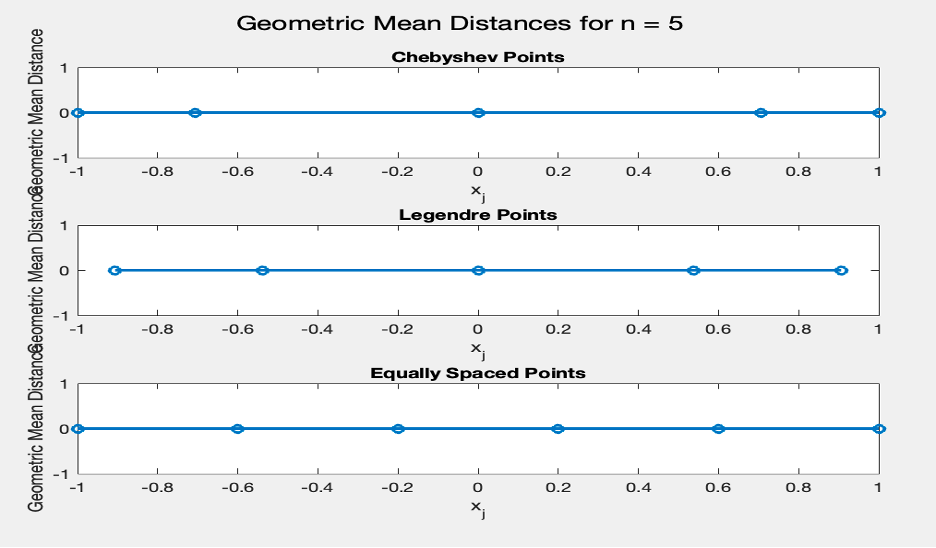}
    \caption{A geometric mean distance plot for Chebyshev points, Legendre points, and equally spaced points for n=5.}
    \label{fig7}
\end{figure}

\begin{figure}[!t]
    \centering
    \includegraphics[width = 8cm, height=10cm]{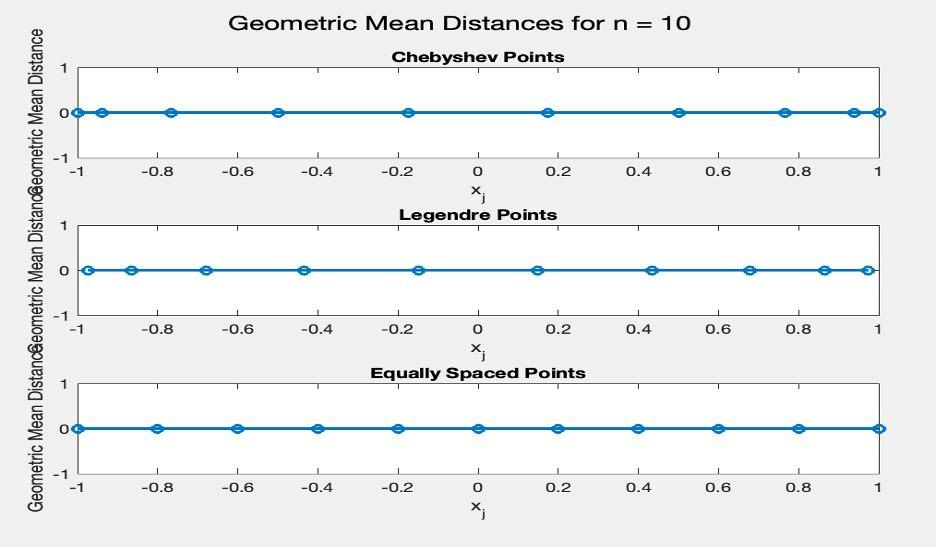}
    \caption{A geometric mean distance plot for Chebyshev points, Legendre points, and equally spaced points for n=10.}
    \label{fig8}
\end{figure}

\begin{figure}[!t]
    \centering
    \includegraphics[width = 8cm, height=10cm]{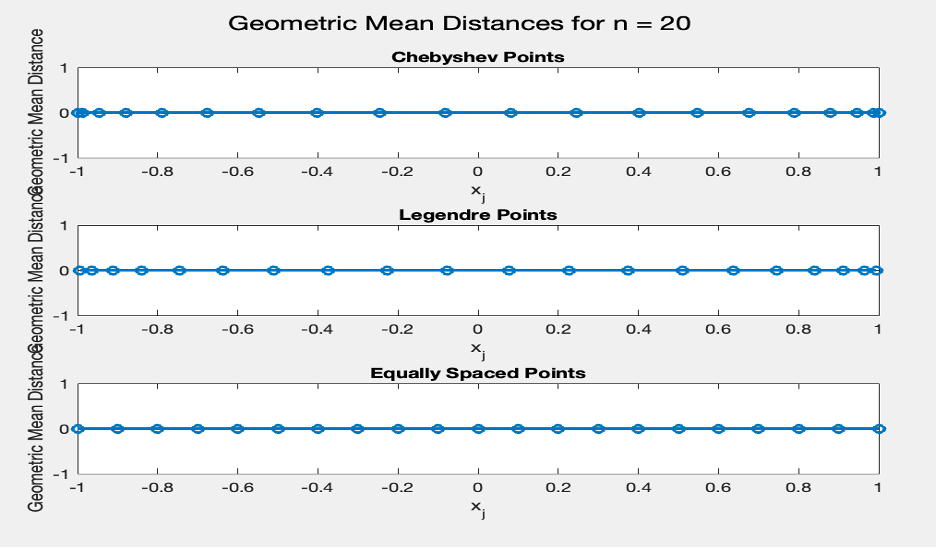}
    \caption{A geometric mean distance plot for Chebyshev points, Legendre points, and equally spaced points for n=20}
    \label{fig9}
\end{figure}

\subsection{Scaled Chebyshev points to the interval [a, b]}
The Chebyshev points are computed in the interval [-1,1] but can be scaled to any interval of preference. Let’s consider Chebyshev points for $n = 9$ in the interval [-1,1] using the function ‘\texttt{chebpts (10)}’ in Matlab. Now using the \texttt{chebfun(@sin,10)} to compute the degree 9 interpolants of $p(x)$ to $sin(x)$ scaled to the interval [-6, 6] and make a semilog plot of $|f(x)-p(x))|$. The results from this exercise will be compared to the same degree 9 interpolants of $p(x)$ to $sin(x)$ but scaled on the interval [0,6].

\begin{figure}[!t]
    \centering
    \includegraphics[width = 8cm, height=10cm]{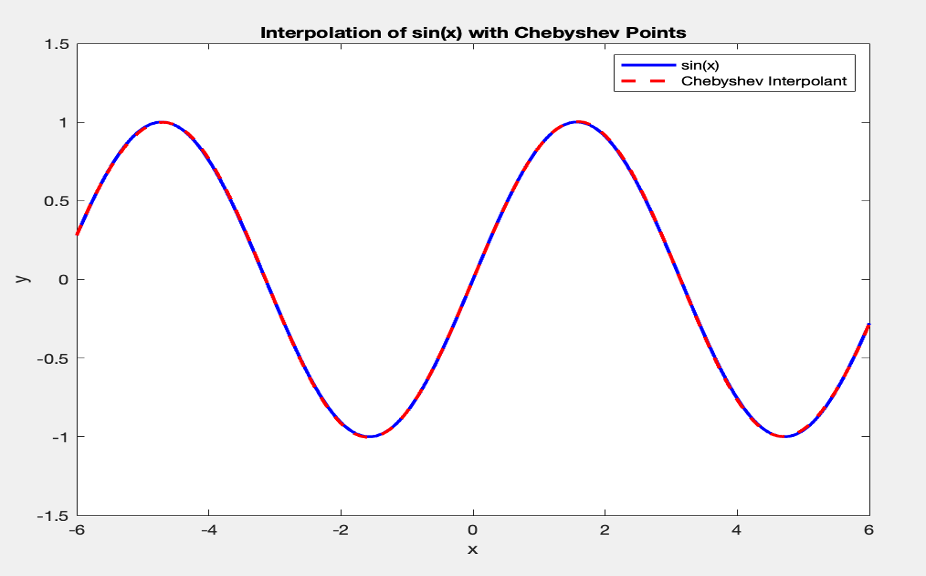}
    \caption{A plot of 9 degrees Chebyshev interpolant of p(x) to sin(x) on the interval [-6, 6]}
    \label{fig10}
\end{figure}

Values of the Chebyshev points for \( n=9 \):
\begin{lstlisting}
    -1.0000
    -0.9397
    -0.766
    -0.500
    -0.1736
     0.1736
     0.500	
     0.7660
     0.9397
     1.00

\end{lstlisting}

Figure \ref{fig10} is the Chebyshev interpolant of $p(x)$ to $sin(x)$. As shown, the Chebyshev points over the interval of [-6, 6] interpolates $sin(x)$ perfectly. However, the scaled interval to [0, 6] plotted over the interval does not (see Figure \ref{fig12}).

\begin{figure}[!t]
    \centering
    \includegraphics[width = 8cm, height=10cm]{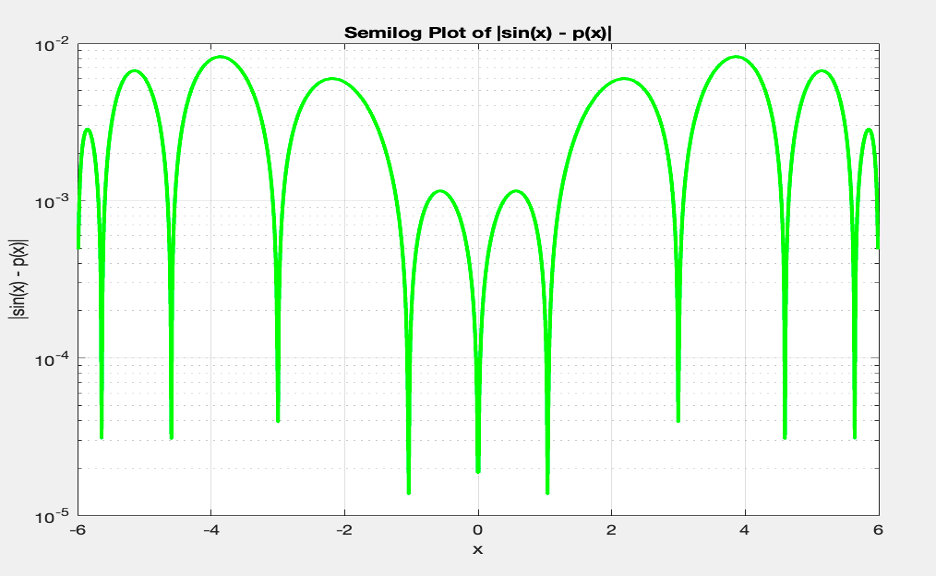}
    \caption{The plot of absolute error of ||f(x)-p(x)||}
    \label{fig11}
\end{figure}

\begin{figure}[!t]
    \centering
    \includegraphics[width = 8cm, height=10cm]{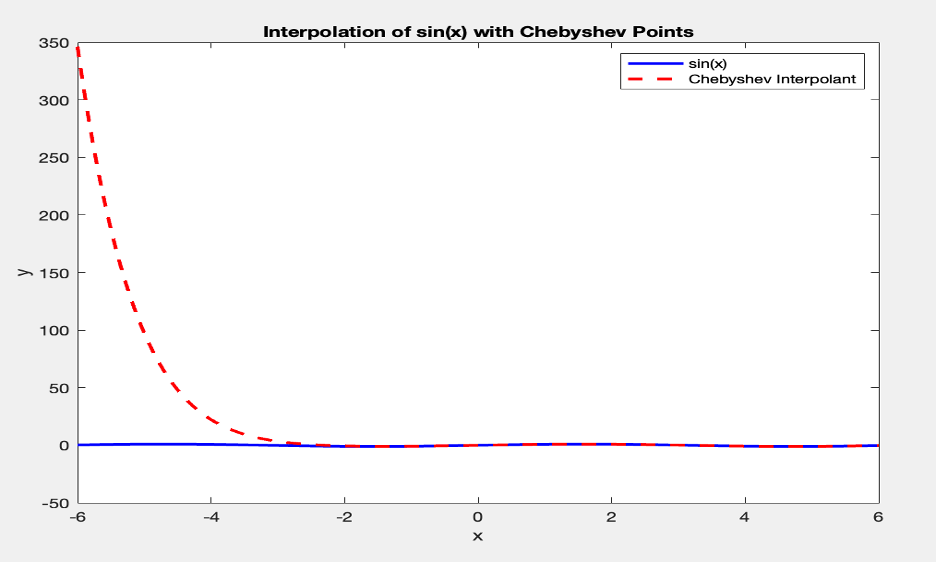}
    \caption{A plot of 9 degrees Chebyshev interpolant of p(x) to sin(x) scaled to the interval [0,6]}
    \label{fig12}
\end{figure}

\begin{figure}[!t]
    \centering
    \includegraphics[width = 8cm, height=10cm]{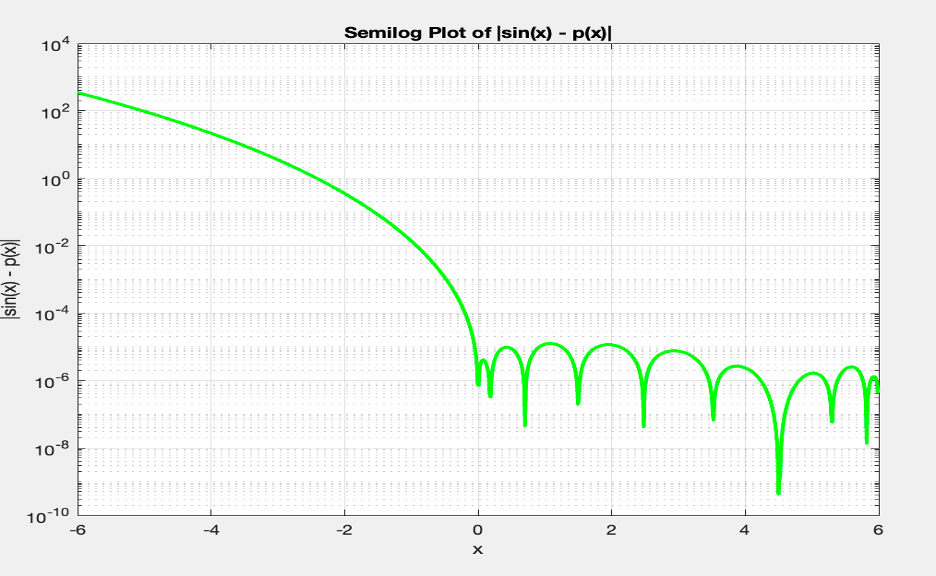}
    \caption{A semilog plot scaled to [0,-6].}
    \label{fig13}
\end{figure}

%%%%%%%%%%%%%%%%%%%%%%%%%%%%%%%%%%%%%%%%%%%%%%%%%%%

\section{Chebyshev polynomials and series}

The interpolation of functions can be done under different fundamental settings and analogies. Some well-known interpolation procedures include those of Fourier and Laurent. In this section, the focus is on the Chebyshev setting with a variable $x$ and a function $f$ defined on [-1, 1].

For
\begin{align}
    x\in\left[-1,1\right],  f\left(x\right)\approx\sum_{k=0}^{n}a_kT_k\left(x\right)
\end{align}

Where $T_k$ is the $kth$ Chebyshev polynomial. The ‘chebpoly(k)’ function in the Chebfun package returns the corresponding Chebfun for $T_k$.

\subsection{A Chebyshev coefficient}

For a function $f(x)=\ \tan^{-1}{(x)}$ on [-1, 1], the coefficient of $T_5$ in the Chebyshev expansion using the chebfun is given below.
\begin{lstlisting}
    x=chebfun('x');
f = atan(x); a = chebpoly(f); format long, a(end:-1:1)'
ans =
0.828427124746190
  -0.047378541243650
   0.004877323527903
  -0.000597726015161
   0.000079763888583
  -0.000011197079759
   0.000001625558989
  -0.000000241714919
   0.000000036592697
  -0.000000005617439
   0.000000000872010
  -0.000000000136603
   0.000000000021562
  -0.000000000003425
   0.000000000000547
  -0.000000000000088
   0.000000000000014
  -0.000000000000002
   0.000000000000000
  -0.000000000000000
\end{lstlisting}
   
Notice that in the chebpoly function, the coefficients for the length of $x$ are given, to call for a specific degree of the Chebeshev coefficients the Matlab codes can be used.

\begin{lstlisting}
    degree = 5; 
chebcoeffs(f, degree + 1)
ans =
   0.828427124746190
  -0.047378541243650
   0.004877323527903
\end{lstlisting}

\subsection{Dependence on wave number}
Let $f(x)$ and $g(x)$ be functions on [-1, 1] defined as 

\begin{align}
    f\left(x\right)= sin (kx)
\end{align}

\begin{align}
    g\left(x\right)=\frac{1}{(+{(kx)}^2)}
\end{align}

\begin{figure}[!t]
    \centering
    \includegraphics[width = 8cm, height=10cm]{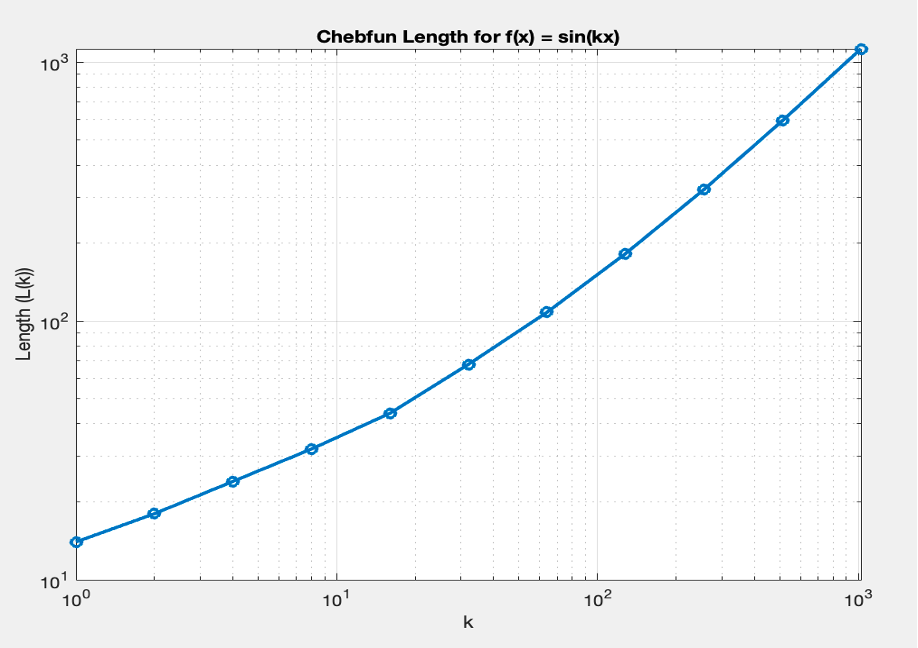}
    \caption{Chebyshev plot of the length $L(k)$ corresponding to $f(x)=sin(kx)$}
    \label{fig14}
\end{figure}

\begin{figure}[!t]
    \centering
    \includegraphics[width = 8cm, height=10cm]{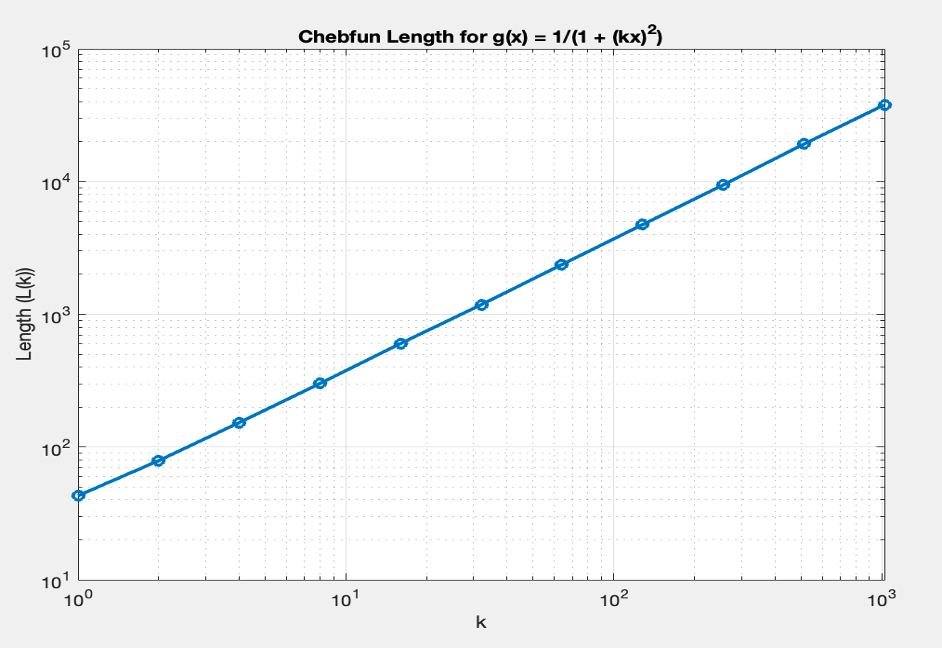}
    \caption{Chebyshev plot of the length L(k) corresponding $g\left(x\right)=\ \frac{1}{(+{(kx)}^2)}$}
    \label{fig15}
\end{figure}

For $k=1,2,4,8…,210$, the length $(L(k))$ of $f(x)$ can be calculated using the chebfun command with the defined function $f(x)$. Figure \ref{fig14} and Figure \ref{fig15} are the log-log plots $L(k)$ as a function of $k$ for $f(x)$ and $g(x)$ respectively. The log-log plots show a linear relationship with the $L(k)$ with the corresponding $f(x)$ showing a slight curve. This gives us valuable insights into how the length of the Chebfun changes with the frequency parameter $k$.

\subsection{Chebyshev series of a complicated function}
Now let's consider the following three functions on [-1, 1].

\begin{align}
  f\left(x\right)=\tanh (x)  
\end{align}

\begin{align}
    g\left(x\right)=\ {10}^{-5}\tanh (10x)
\end{align}

\begin{align}
    h\left(x\right)=\ {10}^{-10}\tanh (100x)
\end{align}

\begin{figure}[!t]
    \centering
    \includegraphics[width = 8cm, height=10cm]{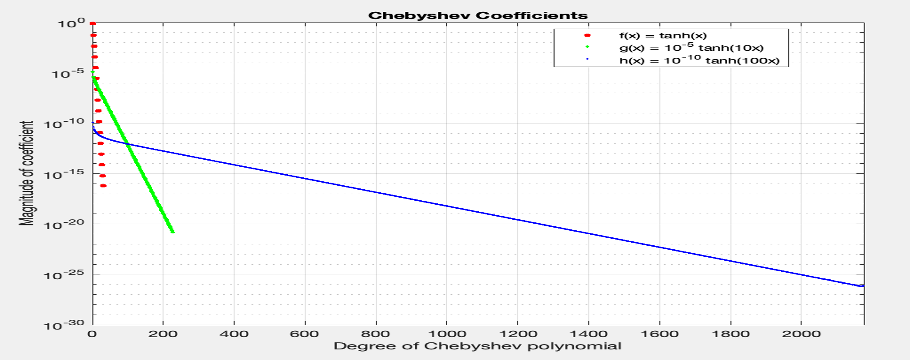}
    \caption{Chebyshev coefficient plot of $f(x), g(x)$, and $h(x)$.}
    \label{fig16}
\end{figure}

\begin{figure}[!t]
    \centering
    \includegraphics[width = 8cm, height=10cm]{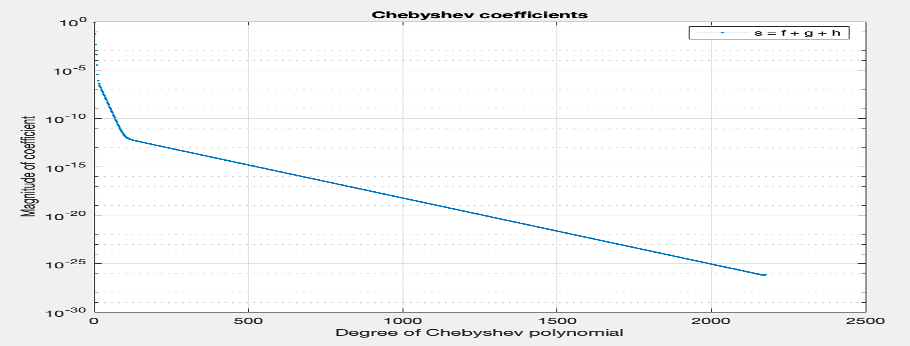}
    \caption{Chebyshev plot for $s = f + g + h$}
    \label{fig17}
\end{figure}

\begin{figure}[!t]
    \centering
    \includegraphics[width = 8cm, height=10cm]{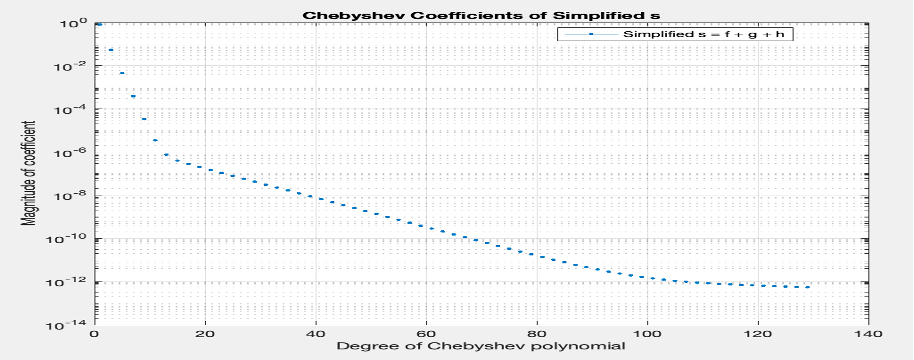}
    \caption{Chebyshev plot of $s = f +g + h$ with the simplify function in chebfun.}
    \label{fig18}
\end{figure}

Making chebfun these functions and using the Chebfun function ‘chebpolyplot’ will produce the results plots in Figure \ref{fig16}. The plots tend to show an interesting pattern as the coefficients of $x$ increase. The function $f(x)$ produces points along the zero points of the degree of Chebyshev polynomial, as the coefficients of $x$ are increased, the functions deviate further from the previous function and become flatter. Figure \ref{fig17} and Figure \ref{fig18} show the plots for $f + g + h$ and when the simplify function in chebfun is applied.

\subsection{Conditioning of the Chebyshev basis}

A quasimatrix, as employed in Chebfun, refers to a structure in the form of a matrix but with a unique property, i.e., one of its dimensions is continuous, while the other remains discrete.
Specifically, a quasimatrix can have more than one column (or, when transposed, more than one row), resulting in an "$\infty $ × $n$" quasimatrix, where each column is a Chebfun. 

\subsubsection{The commands $size(T), cond(T), spy(T)$, and $svd(T)$}

\begin{figure}[!t]
    \centering
    \includegraphics[width = 8cm, height=10cm]{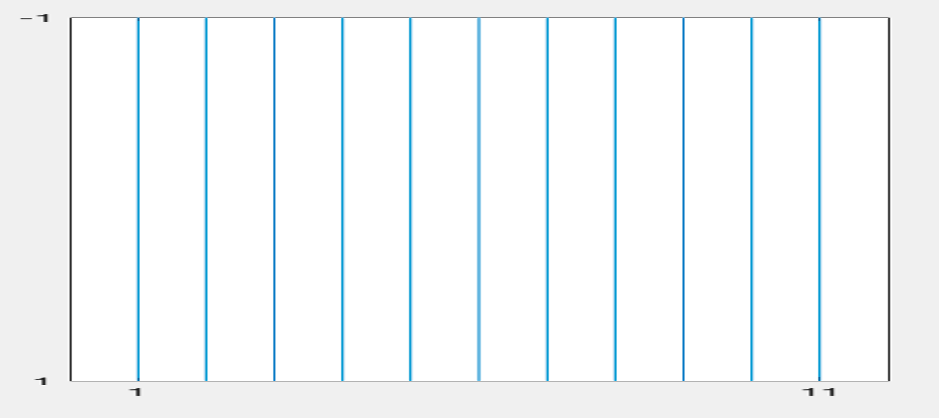}
    \caption{Sparsity plot of T with the $spy(T)$ command}
    \label{fig19}
\end{figure}

Let’s consider a chebfun T (T= chebpoly (0:10)) that is a ‘quasimatrix’. Here, T is an $\infty$ x 11 quasimatrix with each of the 11 columns being a chebfun. The commands size(T), cond(T), spy(T), plot(T), and svd(T) provide information on the size of the quasimatrix, the condition number of the basis, an idea of the shape of a quasimatrix, the plot of the basis function, and the single value decomposition of T.  Figure \ref{fig19} is a plot of the sparsity pattern of the matrix T.

\begin{lstlisting}
T = chebpoly(0:10);
% Size of T (number of columns)
size(T)
ans =
   Inf    11
cond(T)
ans =
   3.712641510134901
svd(T)
ans =
   1.523832995601609
   1.226785409090009
   1.225285973614859
   1.145821182000790
   1.141511078719328
   1.004939431619162
   0.998319547384275
   0.787441811713346
   0.782619113084012
   0.414600480581676
   0.410444421159920
\end{lstlisting}

\begin{figure}[!t]
    \centering
    \includegraphics[width = 8cm, height=10cm]{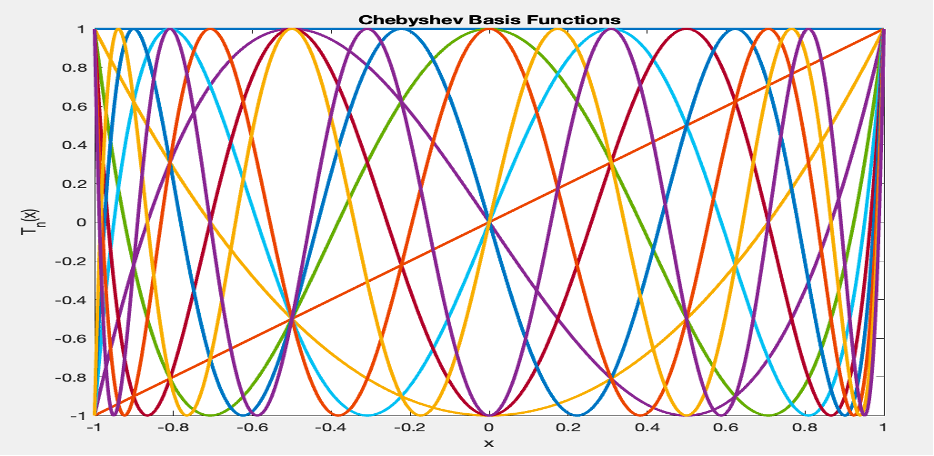}
    \caption{Plot of the Chebyshev basis function of T}
    \label{fig20}
\end{figure}

\subsubsection{Corresponding quasimatrix of monomials}
The corresponding quasimatrix of monomials is executed by the commands in below which resulted in a condition number of 3.072959852624380e+03 for M. 

\begin{lstlisting}
    x = chebfun('x');
M = T;
for j = 0:10
    M(:,j+1) = x.^j;
end
cond(M)
ans =
     3.072959852624380e+03

\end{lstlisting}

\subsubsection{Comparison of the monomial and Chebyshev basis defined on [-1, 1]}

\begin{figure}[!t]
    \centering
    \includegraphics[width = 8cm, height=10cm]{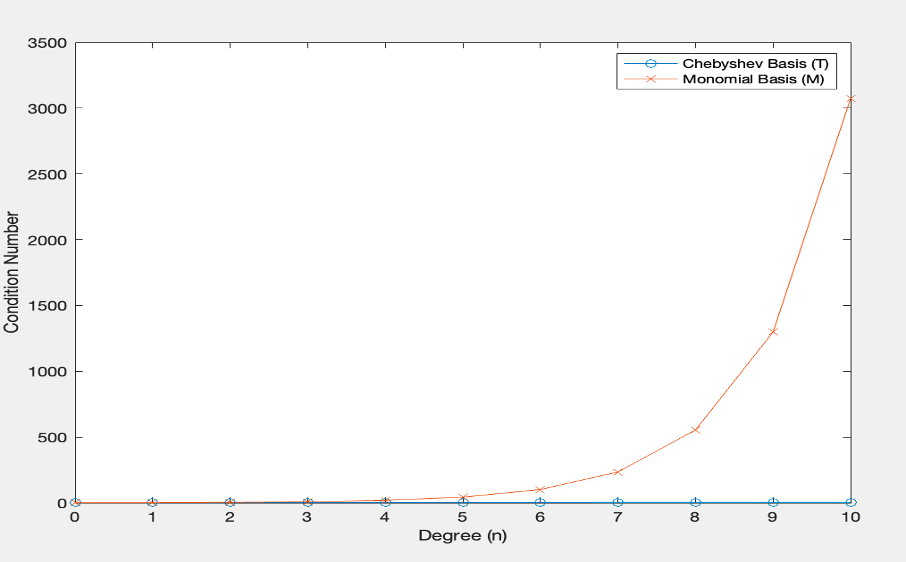}
    \caption{Plot of the monomial and Chebyshev condition numbers}
    \label{fig21}
\end{figure}

The monomial and Chebyshev condition basis is shown in Figure \ref{fig21}. The plot provides valuable insights into the conditioning of the Chebyshev basis (T) and the monomial basis (M) over the range of degrees $n=1, 2, 3, …, 10$. The condition number for the Chebyshev basis remains constant at 0 for $n=1, 2, 3, …, 10$ indicating that the Chebyshev basis is extremely well-conditioned throughout this range. The monomial basis on the other hand starts at 0 for $n=0$ but the condition number starts to rise gradually after $n=3$ to a little over $3000$ indicating an ill conditioning over its range of values. The plot illustrates that the Chebyshev polynomials maintain excellent numerical stability defined on [-1, 1] across different degrees, whereas the monomial basis becomes increasingly prone to numerical instability as the degree of polynomials rises.

\subsubsection{Condition numbers if M is constructed from monomials on [0, 1]}
If M is constructed from the monomials on [0, 1], the condition numbers decreased signifying an improvement in the condition number. This translates to more numerical stability and a less likelihood of being prone to errors.

\begin{lstlisting}
    T=chebpoly(0:10,[0,1]);
x = chebfun('x', [0, 1]);
M = T;
for j = 0:10
    M(:,j+1) = x.^j;
end
cond(M)
ans =  
       2.2871e+07

\end{lstlisting}

\subsection{A function neither even nor odd}

\begin{figure}[!t]
    \centering
    \includegraphics[width = 8cm, height=10cm]{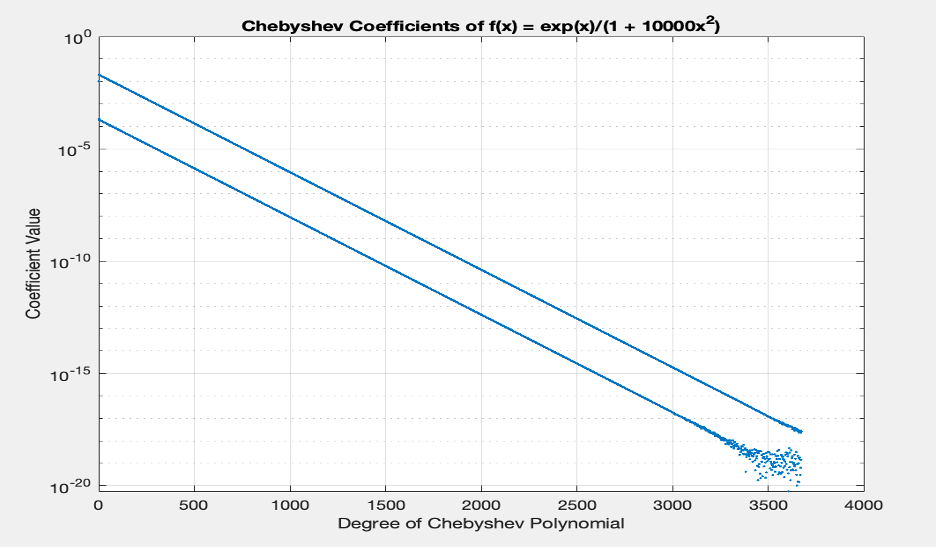}
    \caption{The chebpolyplot of f(x). The appearance of the strips is due to the higher degree oscillations of f(x).}
    \label{fig22}
\end{figure}

Let consider a function \( f(x) = \frac{\exp(x)}{(1+10000x)^2} \) for which we will apply the command \texttt{chebpolyplot} to assess the behavior of the Chebyshev coefficients. In Figure \ref{fig22}, the plot has the appearance of a stripe because of the higher-order degree of oscillation for the Chebyshev coefficients of \( f(x) \).

\subsection{Extrema and roots of Chebyshev polynomials}
The Chebyshev polynomial of the first kind, denoted as $T_{n\left(x\right)}$, has known extrema and roots in the interval [-1, 1]. Here are the formulas for the extrema and roots of $T_{n\left(x\right)}$. The extrema of $T_{n\left(x\right)}$ occur at $n+1$ equally spaced points in the interval [-1, 1]. These extrema are given by:

\begin{align}
    x_k=\cos{\left(\frac{k\pi}{n}\right),\ for\ k=0,\ 1,\ 2,\ \ldots,\ n}
\end{align}

The roots of $T_{n\left(x\right)}$ also occur at n equally spaced points in the interval [-1, 1]. These roots are given by:

\begin{align}
    x_k=\frac{(2k+1)\pi}{2n}\ for\ k=0,\ 1,\ 2,\ \ldots,\ n-1
\end{align}
%%%%%%%%%%%%%%%%%%%%%%%%%%%%%%%%%%%%%%%%%%%%%%%%%%%%%%%%%%%%%%%%%%%%%%

\section{The Chebyshev Interpolation of the Gamma Variate curve}
In our exploration of approximation methods, the robustness of Chebyshev polynomials has stood out in comparison to other techniques. This section specifically delves into the interpolation of the gamma variate function using Chebyshev polynomials, drawing comparisons with Fourier polynomial interpolation. Two distinct approaches are employed for Chebyshev fitting: one with unevenly distributed nodes and another with evenly distributed nodes.

\subsection{Even and unevenly distributed Chebyshev nodes}
Chebyshev nodes, defined by:
\begin{align}
     x_j=cos\left(\frac{(2j+1)\pi}{2n}\right) For j = 0, 1…, n-1
\end{align}

exhibit uneven spacing. However, in terms of $t_k$ defined as
\begin{align}
    t_k=\cos^{-1}{(x_j)=}\frac{(2j+1)\pi}{2n}
\end{align}	

they become evenly spaced.
Consider a gamma variate function with $\alpha = 2$ and $\beta= 1$.

\begin{align}
   h(t) = A(t - t_0)^\alpha \exp\left(-\frac{(t - t_0)}{\beta}\right), \quad t \geq t_0
\end{align}

For the Chebyshev interpolation, we employ the chebfun command on the gamma variate function with evenly distributed nodes $(t = [0,\ \pi,\ 2\pi,\ ...\theta])$. Similarly, for unevenly distributed nodes, t is generated using the command \texttt{sort (rand (1,1000))} * $\pi$. Fourier trigonometry polynomials for even nodes are interpolated using the ‘\texttt{interpft}’ command, a pre-defined function for trigonometric polynomials (see appendix B).

\subsection{Results and Comparison}

\begin{figure}[!t]
    \centering
    \includegraphics[width = 8cm, height=10cm]{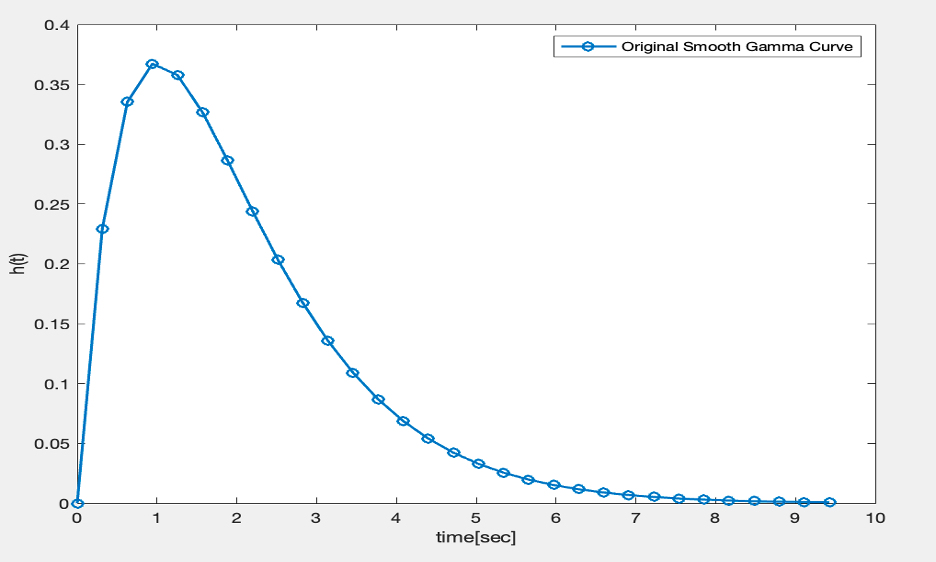}
    \caption{Plot of the original gamma variate function with $\alpha =2$, and $\beta = 1$}
    \label{fig23}
\end{figure}

\begin{figure}[!t]
    \centering
    \includegraphics[width = 8cm, height=10cm]{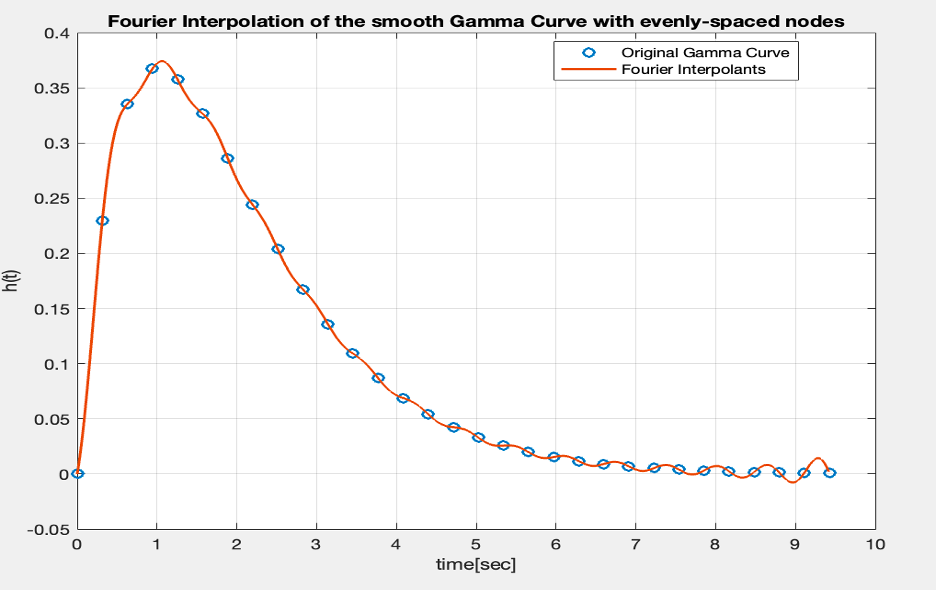}
    \caption{Even sampled Fourier interpolation of the gamma variate function with $\alpha =2$, and $\beta = 1$}
    \label{fig24}
\end{figure}

\begin{figure}[!t]
    \centering
    \includegraphics[width = 8cm, height=10cm]{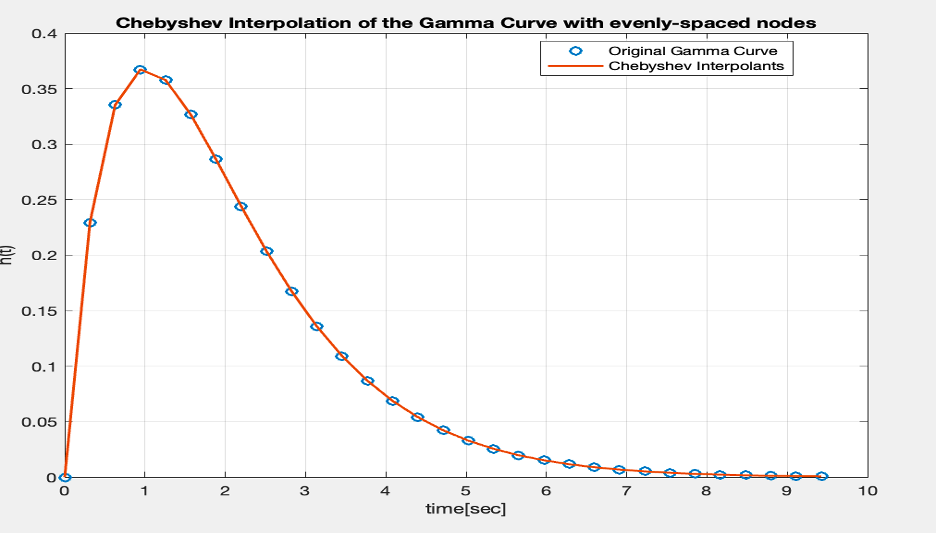}
    \caption{Even sampled Chebyshev interpolation of the gamma variate function with $\alpha =2$, and $\beta = 1$}
    \label{fig25}
\end{figure}

\begin{figure}[!t]
    \centering
    \includegraphics[width = 8cm, height=10cm]{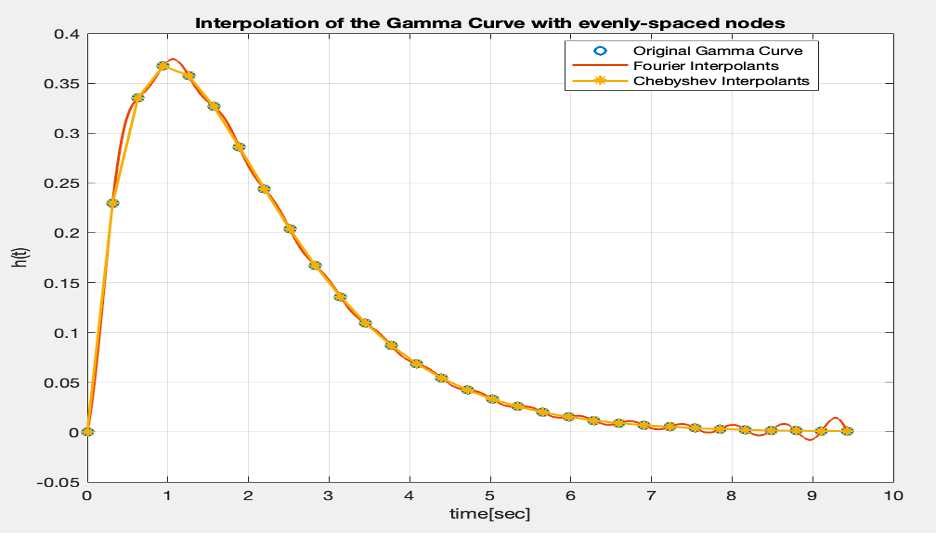}
    \caption{Comparison of the Interpolation of the gamma-variate curve. The maximum value from the Fourier polynomial interpolation and the Chebyshev interpolant for even nodes yields the exact maximum value from the original gamma variate curve}
    \label{fig:26}
\end{figure}

The interpolation outcomes for both even and uneven distributions of Chebyshev nodes, alongside Fourier polynomials, are illustrated in Figures \ref{fig23} to \ref{fig31}. For evenly spaced nodes, both Fourier and Chebyshev interpolations demonstrate excellence. However, in the case of unevenly distributed nodes, only the Chebyshev interpolation is applicable since the Fourier interpolation is only done on the assumption of equally spaced nodes. The unique properties of Chebyshev polynomials make it an ideal interpolation method for both evenly and unevenly distributed nodes.

%%%%%%%%%%%%%%%%%%%%%%%%%%%%%%%%%%%%%%%%%%%%%%%%%%%%%%%%%%%%%%%%%%%%%%%%%

\section{Interpolation of the Gamma Variate curve with noise for evenly and unevenly distributed nodes}

In this section, we focus on the interpolation of the Gamma Variate curve involving noise. The primary objective is to explore and analyze the performance of the Fourier and Chebyshev interpolation in the presence of noise. For noise addition, a Gaussian random noise of 2\% was used.

\subsection{Evenly distributed nodes with noise}

\begin{figure}[!t]
    \centering
    \includegraphics[width = 8cm, height=10cm]{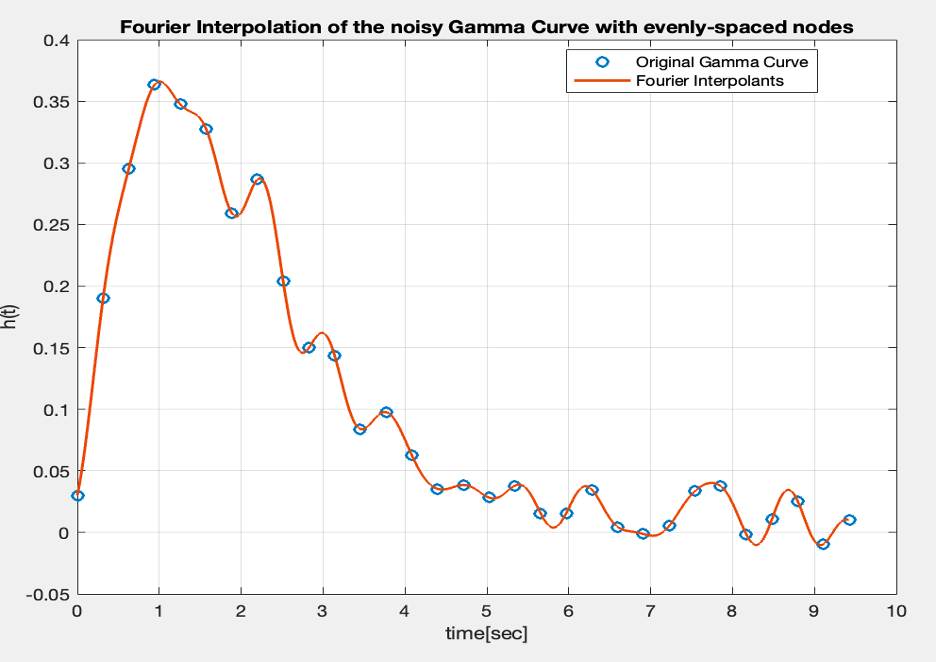}
    \caption{Fourier Interpolation of the gamma-variate curve with added noise for evenly distributed nodes. The original maximum value = 0.3992, Fourier interpolation maximum value = 0.4006.}
    \label{fig27}
\end{figure}

\begin{figure}[!t]
    \centering
    \includegraphics[width = 8cm, height=10cm]{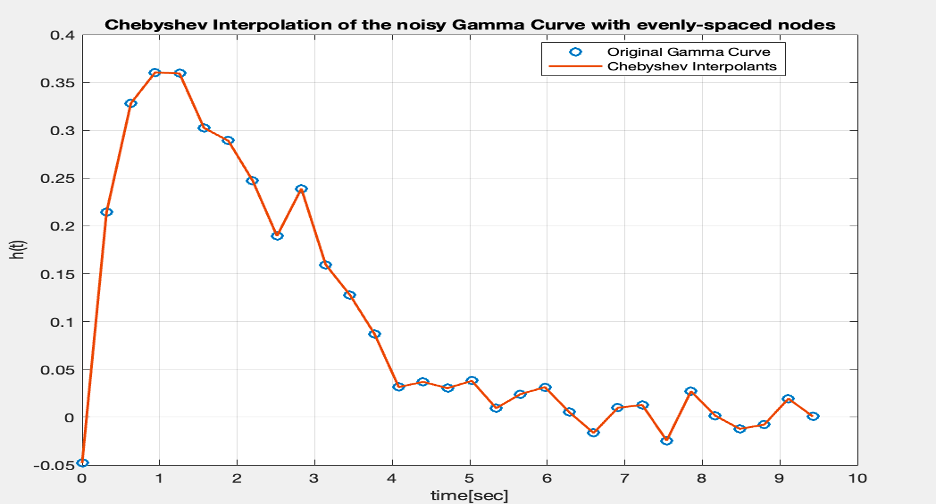}
    \caption{: Chebyshev Interpolation of the gamma-variate curve with added noise for evenly distributed nodes. The original maximum value = 0.3992, Chebyshev interpolation maximum value = 0.3992.}
    \label{fig28}
\end{figure}

\begin{figure}[!t]
    \centering
    \includegraphics[width = 8cm, height=10cm]{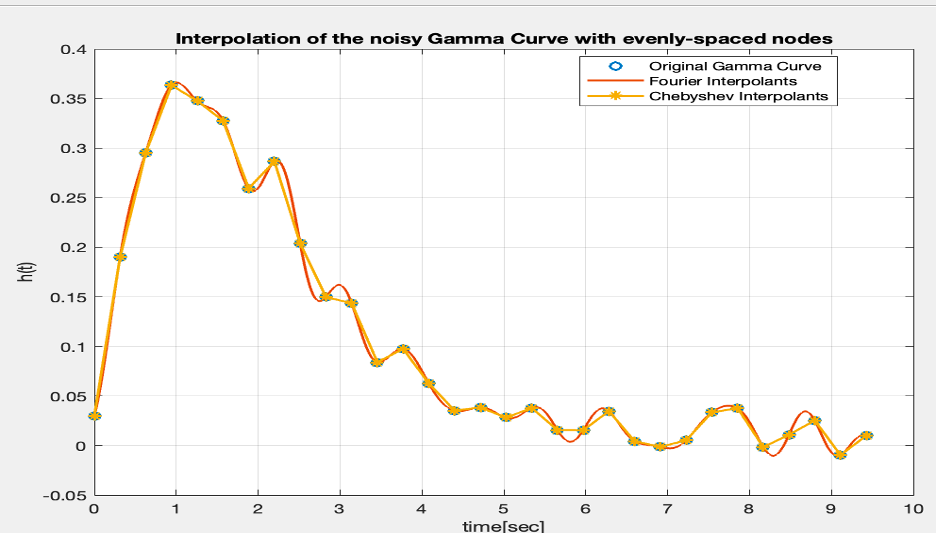}
    \caption{Comparison of the Chebyshev and Fourier Interpolation of the gamma-variate curve with added noise for evenly distributed nodes}
    \label{fig29}
\end{figure}

\begin{figure}[!t]
    \centering
    \includegraphics[width = 8cm, height=10cm]{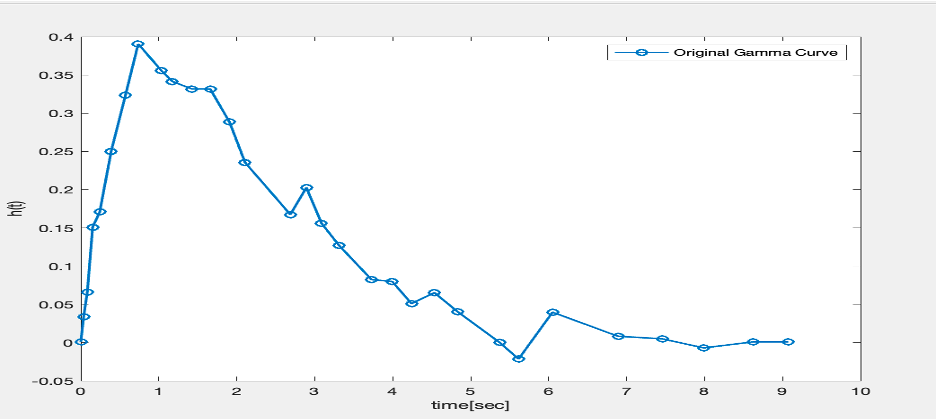}
    \caption{Plot of Noisy a Gamma curve with unevenly distributed nodes for $\alpha=2$, $\beta = 1$}
    \label{fig30}
\end{figure}

\begin{figure}[!t]
    \centering
    \includegraphics[width = 8cm, height=10cm]{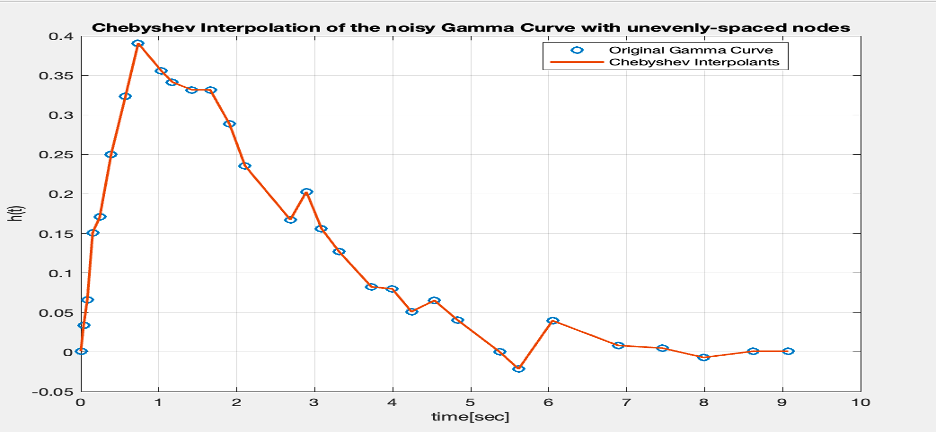}
    \caption{Chebyshev interpolation of a noisy gamma curve with unevenly distributed nodes for $\alpha=2$, $\beta = 1$.}
    \label{fig31}
\end{figure}

For evenly distributed nodes, the Fourier, while largely successful in the interpolation of the gamma curve, exhibits shortcomings by missing the curve's maximum slightly but passes through all data points as expected by the Fourier interpolations. The Chebyshev interpolation on the other hand is excellent. It passes through each point of the Gamma Variate curve with precision and reproduces the original maximum value. Figure \ref{fig27}- \ref{fig29} depicts the comparison of the interpolations to the original gamma curve.

\subsection{Unevenly distributed nodes with noise}

Under unevenly distributed nodes, the Fourier interpolation is not applicable since it violates one of the core assumptions (i.e., equally spaces nodes) necessary for the Fourier interpolation. On the contrary, the Chebyshev interpolation exhibits exceptional accuracy even for unevenly spaced nodes and produces a maximum value equal to that of the original gamma curve. This makes the Chebyshev interpolation a superior choice for interpolation irrespective of the distribution of the nodes as far as they are in the interval [-1, 1]. Figure 30 is the original noisy gamma variate curve while Figure 31 shows the reconstruction using the Chebyshev interpolation. As shown, the Chebyshev interpolation procedure can reconstruct the original curve without any hindrance to how the nodes are sampled with or without noise.

\begin{figure}[!t]
    \centering
    \includegraphics[width = 8cm, height=10cm]{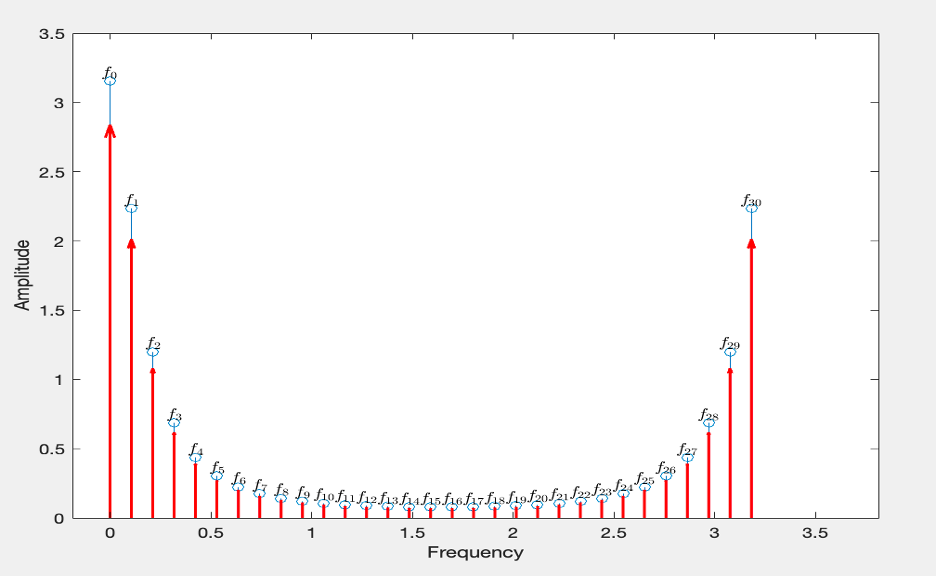}
    \caption{Amplitude plot}
    \label{fig32}
\end{figure}

\begin{figure}[!t]
    \centering
    \includegraphics[width = 8cm, height=10cm]{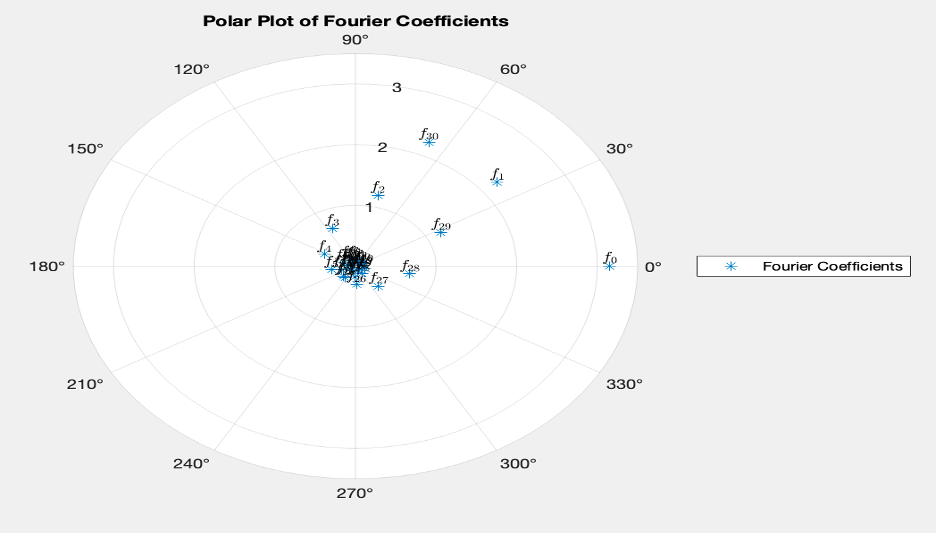}
    \caption{Polar plot of Fourier Coefficients}
    \label{fig33}
\end{figure}

\begin{figure}[!t]
    \centering
    \includegraphics[width = 8cm, height=10cm]{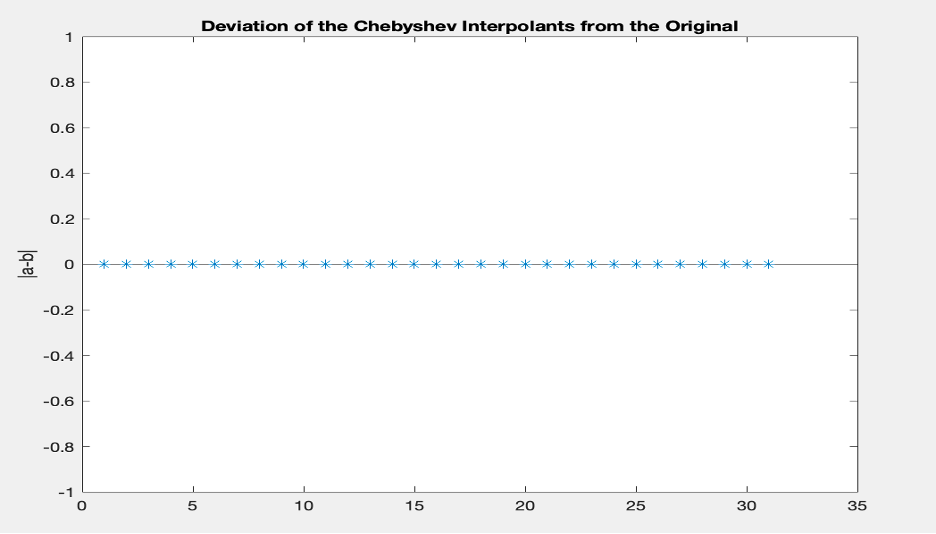}
    \caption{Absolute deviation of the Chebyshev interpolants from the original curve}
    \label{fig34}
\end{figure}

\begin{figure}[!t]
    \centering
    \includegraphics[width = 8cm, height=10cm]{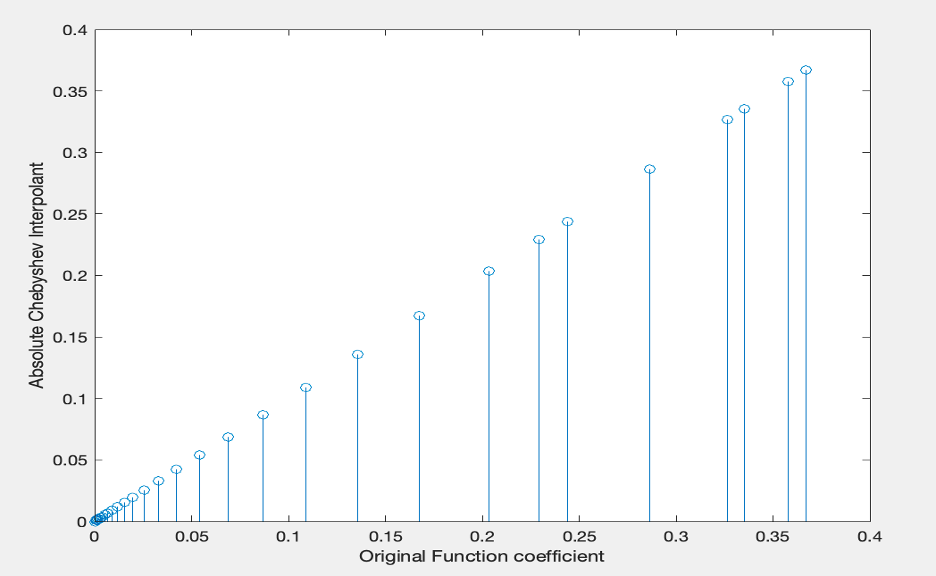}
    \caption{Plot of the absolute Chebyshev interpolants versus the original function}
    \label{fig35}
\end{figure}

\section{Signal Filtering}
\begin{figure}[!t]
    \centering
    \includegraphics[width = 8cm, height=10cm]{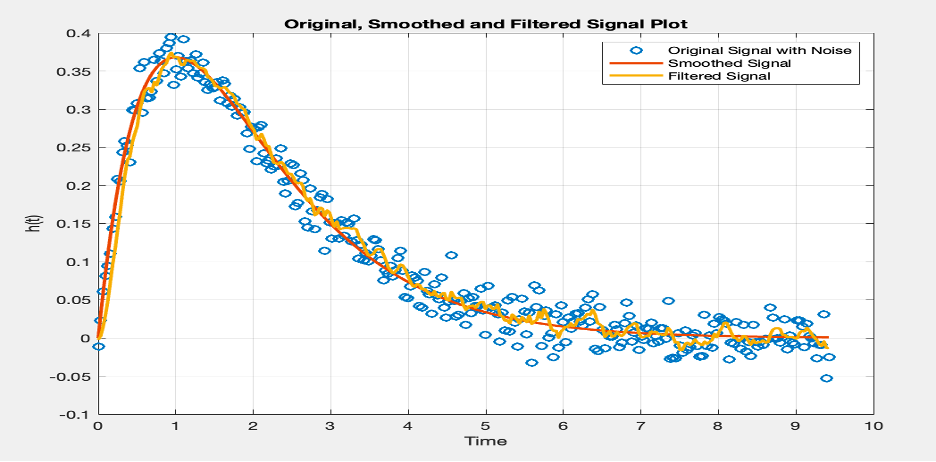}
    \caption{A moving average (window size = 5) filtered signal plot using even nodes.}
    \label{fig36}
\end{figure}

\begin{figure}[!t]
    \centering
    \includegraphics[width = 8cm, height=10cm]{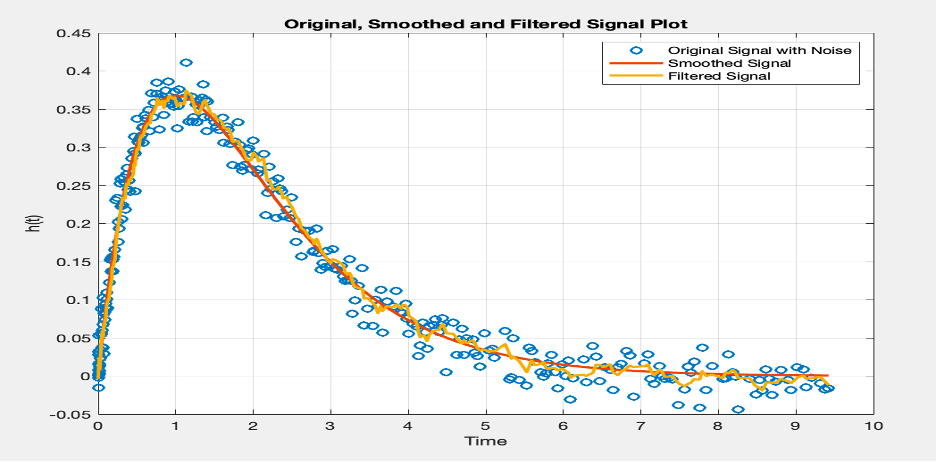}
    \caption{A moving average filtered signal plot using uneven nodes. }
    \label{fig37}
\end{figure}

Noise filtering is a technique performed to remove noise from a set of data. Real-world data is usually characterized by noise, so the filtering of such noise is used to modify or extract key features from such noisy signals. The process adopted here involves the application of a filter to a sequence of data points, where the filter is designed to drop high-frequency coefficients. Figure \ref{fig36} and Figure \ref{fig37} show the filtered gamma variate curves for even and uneven nodes using the ‘\texttt{filter}’ function (a moving average design) in Matlab with a window size of 5.

Filtering is a necessary tool for signal processing, enabling the refinement and extraction of valuable information from signals characterized by noise. However, the selection of an appropriate filter design or window size is key to attaining an intended smooth filtered signal. For this study, different window sizes were considered and the best one closest to the original noise-free curve was selected.

\section{Discussions}
The Chebyshev interpolation has a unique characteristic compared to other polynomial approximations, in that it affords one the comfort to interpolate curves using either of the two parameters $n$ and $\theta$. The interpolation can be expressed in terms of Chebyshev nodes in the $n$-parameter space, where $n$ is the number of data points or nodes used for interpolation for unevenly spaced data points. Alternatively, when the data points are evenly spaced (for angular parameter $\theta$), the Chebyshev interpolation can still be used for such data forms making the Chebyshev interpolation an ideal approximation method for specific characteristics of the dataset at hand. 
Extending the case study by introducing noise further attests to the robustness of the Chebyshev interpolation. Irrespective of how the dataset is formulated (i.e., with even nodes or uneven nodes) or the level of noise, the Chebyshev fitting does not lack precision as in the case of the Fourier fitting which appears to be highly influenced by the uneven distribution of the nodes. This is due to the violation of one of the basic requirements for the Fourier interpolations to be successful, which is the need for uniformly distributed nodes. However, in real-world scenarios, many datasets involve unevenly distributed nodes. The Chebyshev interpolations eliminate this limitation of the Fourier interpolations and maintain a high level of accuracy and efficiency.  Though the Chebyshev interpolation is efficient for any node distribution, its major limitation is the need for the function to be clustered between $\pm\ 1$. 

Considering the above advantages that the Chebyshev interpolation has over the Fourier interpolation, it is worth noting that the Fourier interpolation is widely used over the Chebyshev interpolation in perfusion analysis. Perhaps because the Fourier method offers alternative methods via convolution and deconvolution to reconstruct the true curve of a function amid noise.

As displayed in Figure \ref{fig34} and Figure \ref{fig35}, the deviation of the Chebyshev interpolants from the original function shows a perfect estimation of the original curve with a mean absolute deviation of zero.

\section{Conclusion}
 The Chebyshev polynomials interpolate evenly-spaced and unevenly-spaced points perfectly, particularly with clustering around $-1$ and $1$. However, when scaling is done in $[0, 6]$, the interpolation isn’t as perfect. A comparison of the Chebyshev points with the Legendre points reveals significant similarities, with each point having approximately the same average distance from other points in the geometric mean setting. Again, the Chebyshev polynomial is robust and well-conditioned over its range of values.
Chebyshev interpolation of the gamma variate function proves highly accurate, even when the case of unevenly distributed nodes was tried making it a reliable interpolation method for different node structures compared to Fourier polynomial interpolation which is only usable for equally spaced nodes.
The Chebyshev interpolants provide an excellent approximation for both evenly-spaced nodes and unevenly-spaced nodes making it an ideal approximation method for numerical analysis even amidst noise.
\section{Future work}
This research centered on utilizing Chebyshev interpolation, with a specific emphasis on employing the chebfun function as an alternative for interpolating perfusion data. Notably, the study demonstrated the excellent performance of Chebyshev polynomials in interpolating the gamma curve, a renowned function for tracer dilutions for two scenarios, i.e., for evenly spaced and unevenly spaced nodes with and without noise.

Future investigations should expand their scope by reconstructing the true underlying Gamma curve amid noise with the Chebyshev interpolants and the results compared to the reconstruction by the Fast Fourier Transform (FFT).

\bibliographystyle{ieeetr}
\bibliography{mybibfile}
\nocite*

%%%%%%%%%%%%%%%%%%%%%%%%%%%%%%%%%%%%%%%%%%%%%%%%%%%%%%%%%%%%%%%%%%%%%%%

\section*{Appendix A}
\section*{MATLAB codes}
\begin{lstlisting}
% Figure 1 
tic, p = chebfun(2*rand(10,1)-1), toc; 
figure;
hold off, plot(p,'-'), hold on, plot(p,'.k')
ylim([-1.7 1.7]), grid on
title('Chebyshev interpolant through random data', 'fontsize', 9)
figure;
plot(p, '*-', 'interval', [0.9999 1], 'linewidth', 3)
\end{lstlisting}

\begin{lstlisting}
% Figure 2 
tic, p = chebfun(2*rand(100,1)-1), toc;
figure;
hold off, plot(p, '-'), hold on, plot(p, '.k')
ylim([-1.7 1.7]), grid on
title('Chebyshev interpolant through random data', 'fontsize', 9)
plot(p, '*-', 'interval', [0.9999 1], 'linewidth', 3)
\end{lstlisting}

\begin{lstlisting}
% Figure 3 
tic, p = chebfun(2*rand(1000,1)-1), toc;
hold off, plot(p, '-'), hold on, plot(p, '.k')
ylim([-1.7 1.7]), grid on
title('Chebyshev interpolant through random data', 'fontsize', 9)
figure;
plot(p, '.-', 'interval', [0.9999 1], 'linewidth', 3)
\end{lstlisting}

\begin{lstlisting}
% Figure 4
tic, p = chebfun(2*rand(10000,1)-1), toc;
hold off, plot(p, '-'), hold on, plot(p, '.k')
ylim([-1.7 1.7]), grid on
title('Chebyshev interpolant through random data', 'fontsize', 9)
plot(p, '*-', 'interval', [0.9999 1], 'linewidth', 3)
\end{lstlisting}

\begin{lstlisting}
% Figure 5
n = 99;
p1 = chebpts(n+1);
p2 = legpts(n+1);
pts_diff = abs(p1 - p2);
max_diff = max(pts_diff);
figure;
plot(p1, 'o');
hold on;
plot(p2, '-', 'lineWidth', 1.5);
legend('Chebyshev Points of the First Kind', 'Legendre Points', 'Location', 'best');
xlabel('n');
hold off;
title('Comparison of the Chebyshev Points of the First Kind and the Legendre Points');
max_diff
\end{lstlisting}

\begin{lstlisting}
% Figure 6
ni = 1:100;
fx1 = 'exp(x)';
fx2 = '1/(1+25*x.^2)';
f1 = chebfun(fx1);
f2 = chebfun(fx2);
for n = ni
    e(n) = norm(f1 - chebfun(fx1, n));
end
figure;
subplot(2, 1, 1);
semilogy(e), title('semilog scale for f(x) = exp(x)')
xlabel('n')
hold on;
for n = ni
    e(n) = norm(f2 - chebfun(fx2, n));
end
subplot(2, 1, 2);
semilogy(e), title('semilog scale for f(x) = 1/(1+25*x.^2)')
xlabel('n')
hold off;
\end{lstlisting}

\begin{lstlisting}
% Finding the machine precision value
machine_epsilon = eps; % Get the machine epsilon value
n_max = 1000;
n = 1;
while n <= n_max
    e1 = norm(f1 - chebfun(fx1, n));
    e2 = norm(f2 - chebfun(fx2, n));
    if e1 < machine_epsilon && e2 < machine_epsilon
        fprintf('For machine precision, n must be at least %d\n', n);
        break;
    end
    n = n + 1;
end
if n > n_max
    fprintf('Machine precision not reached within n_max iterations.\n');
end
\end{lstlisting}

\begin{lstlisting}
% Figure 8
meandistance(10);
\end{lstlisting}

\begin{lstlisting}
% Figure 9
meandistance(20);
\end{lstlisting}

\begin{lstlisting}
% Figure 10
% Compute Chebyshev interpolant to sin(x)
f = chebfun(@sin, [-6, 6], 10);

% Plot sin(x) and interpolant p(x) over the larger interval [-6, 6]
x_vals = linspace(-6, 6, 1000);
figure;
plot(x_vals, sin(x_vals), 'b', x_vals, f(x_vals), 'r--', 'LineWidth', 2);
legend('sin(x)', 'Chebyshev Interpolant');
title('Interpolation of sin(x) with Chebyshev Points');
xlabel('x');
ylabel('y');
\end{lstlisting}

\begin{lstlisting}
% Figure 11
% Semilog plot of absolute error
figure;
semilogy(x_vals, abs(sin(x_vals) - f(x_vals)), 'g', 'LineWidth', 2);
title('Semilog Plot of |sin(x) - p(x)|');
xlabel('x');
ylabel('|sin(x) - p(x)|');
grid on;
\end{lstlisting}

\begin{lstlisting}
% Figure 12
% Compute Chebyshev interpolant to sin(x)
f = chebfun(@sin, [0, 6], 10);

% Plot sin(x) and interpolant p(x) over the interval [0, 6]
x_vals = linspace(-6, 6, 1000);
figure;
plot(x_vals, sin(x_vals), 'b', x_vals, f(x_vals), 'r--', 'LineWidth', 2);
legend('sin(x)', 'Chebyshev Interpolant');
title('Interpolation of sin(x) with Chebyshev Points');
xlabel('x');
ylabel('y');
\end{lstlisting}

\begin{lstlisting}
% Figure 13
% Semilog plot of absolute error
figure;
semilogy(x_vals, abs(sin(x_vals) - f(x_vals)), 'g', 'LineWidth', 2);
title('Semilog Plot of |sin(x) - p(x)|');
xlabel('x');
ylabel('|sin(x) - p(x)|');
grid on;
\end{lstlisting}

\begin{lstlisting}
% Figure 14
% Define the function for f(x) = sin(kx)
f = @(k, x) sin(k * x);
% Values of k
k_values = 2.^(0:10);
% Calculate the Chebfun lengths for each k
L_values_f = zeros(size(k_values));
for i = 1:length(k_values)
    k = k_values(i);
    f_k = @(x) f(k, x);
    F_k = chebfun(f_k, [-1, 1]);
    L_values_f(i) = length(F_k);
end

% Log-log plot
figure;
loglog(k_values, L_values_f, 'o-', 'LineWidth', 2);
title('Chebfun Length for f(x) = sin(kx)');
xlabel('k');
ylabel('Length (L(k))');
grid on;
\end{lstlisting}

\begin{lstlisting}
% Figure 15
% Define the function for g(x) = 1/(1 + (kx)^2)
g = @(k, x) 1./(1 + (k * x).^2);

% Calculate the Chebfun lengths for each k
L_values_g = zeros(size(k_values));
for i = 1:length(k_values)
    k = k_values(i);
    g_k = @(x) g(k, x);
    G_k = chebfun(g_k, [-1, 1]);
    L_values_g(i) = length(G_k);
end

% Log-log plot
figure;
loglog(k_values, L_values_g, 'o-', 'LineWidth', 2);
title('Chebfun Length for g(x) = 1/(1 + (kx)^2)');
xlabel('k');
ylabel('Length (L(k))');
grid on;
\end{lstlisting}

\begin{lstlisting}
% Figure 16
f = chebfun(@tanh, [-1, 1]);
g = 1e-5 * chebfun(@(x) tanh(10 * x), [-1, 1]);
h = 1e-10 * chebfun(@(x) tanh(100 * x), [-1, 1]);
% Create Chebpoly plots for the coefficients
figure;
chebpolyplot(f, '.-', 'DisplayName', 'f(x) = tanh(x)');
hold on;
chebpolyplot(g, '.-', 'DisplayName', 'g(x) = 10^{-5} tanh(10x)');
chebpolyplot(h, '.-', 'DisplayName', 'h(x) = 10^{-10} tanh(100x)');
title('Chebyshev Coefficients');
legend('Location', 'Best');
\end{lstlisting}

\begin{lstlisting}
% Figure 17
% Define the function s as the sum of f, g, and h
s = f + g + h;
% Create a Chebpoly plot for the coefficients of s
figure;
chebpolyplot(s, '.-', 'DisplayName', 's = f + g + h');
title('Chebyshev Coefficients of s');
legend('Location', 'Best');
\end{lstlisting}

\begin{lstlisting}
% Figure 18
% Apply the simplify command to s to remove the tail
simplified_s = simplify(s);

% Create a Chebpoly plot for the coefficients of simplified s
figure;
chebpolyplot(simplified_s, '.-', 'DisplayName', 'Simplified s = f + g + h');
title('Chebyshev Coefficients of Simplified s');
legend('Location', 'Best');
\end{lstlisting}

\begin{lstlisting}
% Figure 19
T = chebpoly(0:10);
size_T = size(T);
figure;
spy(T);

% The following plots are related to the sparse structure of T.
\end{lstlisting}

\begin{lstlisting}
% Figure 20
figure;
plot(T);
\end{lstlisting}

\begin{lstlisting}
% Figure 21
singular_values_T = svd(T);
condition_number_T = cond(T);
x = chebfun('x');
M = T;
for j = 0:10
    M(:, j+1) = x.^j;
end
condition_number_M = cond(M);
condition_numbers_T = zeros(11, 1);
condition_numbers_M = zeros(11, 1);

for n = 0:10
    T_n = chebpoly(0:n);
    M_n = T_n;
    for j = 0:n
        M_n(:, j+1) = x.^j;
    end
    condition_numbers_T(n+1) = cond(T_n);
    condition_numbers_M(n+1) = cond(M_n);
end
figure;
plot(0:10, condition_numbers_T, '-o', 0:10, condition_numbers_M, '-x');
xlabel('Degree (n)');
ylabel('Condition Number');
legend('Chebyshev Basis (T)', 'Monomial Basis (M)');
\end{lstlisting}

\begin{lstlisting}
% Figure 22
figure;
f = chebfun(@(x) exp(x) ./ (1 + 10000 * x.^2), [-1, 1]);
chebpolyplot(f, 'DisplayName', 'f(x) = exp(x)/(1 + 10000x^2)');
title('Absolute values of Chebyshev coefficients');
legend('Location', 'Best');
grid on;
\end{lstlisting}

\begin{lstlisting}
% Figure 23
figure;
dx = 3*pi/30;
x = 0:dx:3*pi;
alpha = 2;
beta = 1;
f = (beta^alpha * x.^(alpha - 1) .* exp(-beta * x)) / gamma(alpha);
plot(x, f, '-o', 'DisplayName', 'Original Gamma Curve', 'lineWidth', 1.5);
legend('Location', 'best');
\end{lstlisting}

\begin{lstlisting}
% Figure 24
figure;
plot(x, f, '-o', 'DisplayName', 'Original Gamma Curve', 'lineWidth', 1.5);
hold on
N = 1000;
y = interpft(f, N);
dy = dx * length(x) / N;
x2 = 0:dy:3*pi;
y = y(1:length(x2));
plot(x2, y, '-', 'DisplayName', 'Fourier Interpolants', 'lineWidth', 1.5)
legend('Location', 'best');
xlabel('time[sec]');
ylabel('h(t)');
title('Fourier Interpolation of the smooth Gamma Curve with evenly-spaced nodes');
\end{lstlisting}

\begin{lstlisting}
% Figure 25
figure;
plot(x, f, '-o', 'DisplayName', 'Original Gamma Curve', 'lineWidth', 1.5);
hold on
cf = chebfun(f, [0, 1]);
a = chebpoly(cf);
plot(x, a, '-', 'DisplayName', 'Chebyshev Interpolants', 'lineWidth', 1.5);
legend('Location', 'best');
xlabel('time[sec]');
ylabel('h(t)');
title('Chebyshev Interpolation of the smooth Gamma Curve with evenly-spaced nodes');
grid on;
hold off;
\end{lstlisting}

\begin{lstlisting}
% Figure 26
figure;
plot(x, f, 'o', 'DisplayName', 'Original Gamma Curve', 'lineWidth', 1.5);
hold on
plot(x2, y, '-', 'DisplayName', 'Fourier Interpolants', 'lineWidth', 1.5)
plot(x, a, '*-', 'DisplayName', 'Chebyshev Interpolants', 'lineWidth', 1.5);
legend('Location', 'best');
xlabel('time[sec]');
ylabel('h(t)');
title('Interpolation of the smooth Gamma Curve with evenly-spaced nodes');
grid on;
hold off;
\end{lstlisting}

\begin{lstlisting}
% Figure 27
figure;
dx = 3*pi/30;
x = 0:dx:3*pi;
alpha = 2;
beta = 1;
noise = 0.02 * randn(1, length(x));
f_noise = f + noise;
plot(x, f_noise, 'o', 'DisplayName', 'Original Gamma Curve', 'lineWidth', 1.5);
hold on
N = 1000;
y = interpft(f_noise, N);
dy = dx * length(x) / N;
x2 = 0:dy:3*pi;
y = y(1:length(x2));
plot(x2, y, '-', 'DisplayName', 'Fourier Interpolants', 'lineWidth', 1.5)
legend('Location', 'best');
xlabel('time[sec]');
ylabel('h(t)');
title('Interpolation of the noisy Gamma Curve with evenly-spaced nodes');
grid on;
hold off;
\end{lstlisting}

\begin{lstlisting}
% Figure 28
figure;
plot(x, f_noise, 'o', 'DisplayName', 'Original Gamma Curve', 'lineWidth', 1.5);
hold on
cf = chebfun(f_noise, [0, 1]);
a = chebpoly(cf);
plot(x, a, '-', 'DisplayName', 'Chebyshev Interpolants', 'lineWidth', 1.5);
legend('Location', 'best');
xlabel('time[sec]');
ylabel('h(t)');
title('Chebyshev Interpolation of the noisy Gamma Curve with evenly-spaced nodes');
grid on;
hold off;
\end{lstlisting}

\begin{lstlisting}
% Figure 29
figure;
plot(x, f_noise, 'o', 'DisplayName', 'Original Gamma Curve', 'lineWidth', 1.5);
hold on
plot(x2, y, '-', 'DisplayName', 'Fourier Interpolants', 'lineWidth', 1.5)
plot(x, a, '*-', 'DisplayName', 'Chebyshev Interpolants', 'lineWidth', 1.5);
legend('Location', 'best');
xlabel('time[sec]');
ylabel('h(t)');
title('Interpolation of the noisy Gamma Curve with evenly-spaced nodes');
grid on;
hold off;
\end{lstlisting}

\begin{lstlisting}
% Figure 30
% Unevenly spaced with noise
figure;
dx = 3 * pi / 30;
x = 0:dx:3 * pi;
uneven_x = sort(rand(1, length(x)));
x = x .* uneven_x;
alpha = 2;
beta = 1;
noise = 0.02 * randn(1, length(x));
f = (beta^alpha * x.^(alpha - 1) .* exp(-beta * x)) / gamma(alpha);
f_noise = f + noise;
plot(x, f_noise, '-o', 'DisplayName', 'Original Gamma Curve', 'lineWidth', 1.5);
legend('Location', 'best');
xlabel('time[sec]');
ylabel('h(t)');
title('Original noisy Gamma Curve with unevenly-spaced nodes');
grid on;
\end{lstlisting}

\begin{lstlisting}
% Figure 31
cf = chebfun(f_noise, [0, 1]);
a = chebpoly(cf);
figure;
plot(x, f_noise, '-o', 'DisplayName', 'Original Gamma Curve', 'lineWidth', 1.5);
hold on
plot(x, a, '-', 'DisplayName', 'Chebyshev Interpolants', 'lineWidth', 1.5);
legend('Location', 'best');
xlabel('time[sec]');
ylabel('h(t)');
title('Chebyshev Interpolation of the noisy Gamma Curve with unevenly-spaced nodes');
grid on;
hold off;
\end{lstlisting}

\begin{lstlisting}
% Figure 32
% Amplitude plots
dx = 3 * pi / 30;
x = 0:dx:3 * pi;
alpha = 2;
beta = 1;
f = (beta^alpha * x.^(alpha - 1) .* exp(-beta * x)) / gamma(alpha);
% Perform Fourier interpolation
N = length(x);
y = interpft(f, N);
dy = dx * length(x) / N;
x2 = 0:dy:3 * pi;
y = y(1:length(x2));
% Calculate Fourier coefficients
amplitude = abs(fft(f));
frequencies = linspace(0, 1 / dy, length(amplitude));
% Construct the amplitude plot with coefficients labeled
figure;
stem(frequencies, amplitude, 'DisplayName', 'Fourier Coefficients');
xlabel('Frequency');
ylabel('Amplitude');
\end{lstlisting}

\begin{lstlisting}
% Figure 33
% Polar plots
dx = 3 * pi / 30;
x = 0:dx:3 * pi;
alpha = 2;
beta = 1;
f = (beta^alpha * x.^(alpha - 1) .* exp(-beta * x)) / gamma(alpha);
% Perform Fourier interpolation
N = length(x);
y = interpft(f, N);
dy = dx * length(x) / N;
x2 = 0:dy:3 * pi;
y = y(1:length(x2));
% Calculate Fourier coefficients
amplitude = abs(fft(y));
frequencies = linspace(0, 1 / dy, length(amplitude));
% Create polar plot with coefficients labeled
theta = 2 * pi * frequencies;
rho = amplitude;
figure;
polarplot(theta, rho, '*', 'DisplayName', 'Fourier Coefficients');
title('Polar Plot of Fourier Coefficients');
\end{lstlisting}

\begin{lstlisting}
% Figure 34
figure;
stem(abs(a - f'), '*')
ylabel('|a - b|')
title('Deviation of the Chebyshev Interpolants from the Original')
\end{lstlisting}

\begin{lstlisting}
% Figure 35
% Deviation vs original curve plot
cf = chebfun(f, [0, 1]);
a = chebpoly(cf);
figure;
stem(abs(a), f)
xlabel('Original Function coefficient')
ylabel('Absolute Chebyshev Interpolant')
\end{lstlisting}

\begin{lstlisting}
% Figure 36 and 37
% Filter
% Given parameters
dx = 3 * pi / 300;
x = 0:dx:3 * pi;
alpha = 2;
beta = 1;
noise = 0.02 * randn(1, length(x));
f = (beta^alpha * x.^(alpha - 1) .* exp(-beta * x)) / gamma(alpha);
signal = f + noise;

% Create a sample signal
% Apply a moving average filter for smoothing
window_size = 5;  % size of the moving average window
smoothed_signal = filter(ones(1, window_size) / window_size, 1, signal);

% Plot the original and smoothed signals
figure;
plot(x, signal, 'o', 'LineWidth', 1.5);
hold on;
plot(x, f, '-', 'LineWidth', 2);
plot(x, smoothed_signal, 'LineWidth', 2);
title('Original, Smoothed and Filtered Signal Plot');
legend('Original Signal with Noise', 'Smoothed Signal', 'Filtered Signal');
xlabel('Time');
ylabel('h(t)');
grid on;
hold off;
\end{lstlisting}

\section*{Appendix B}
\section*{MATLAB Functions}
\begin{lstlisting}

%%Meandistance function

function meandistance(n)
    % Function to calculate and plot geometric mean distances for different point sets
    
    % Chebyshev points
    cheb_points = chebpts(n);

    % Legendre points
    legendre_points = legpts(n);

    % Equally spaced points
    equally_spaced_points = linspace(-1, 1, n+1);

    % Function to calculate geometric mean distance for a given point
    calculate_distance = @(x, points) prod(sqrt((x - points).^2));

    % Calculate geometric mean distances for each set of points
    cheb_distances = arrayfun(@(x) calculate_distance(x, cheb_points), cheb_points);
    legendre_distances = arrayfun(@(x) calculate_distance(x, legendre_points), legendre_points);
    equally_spaced_distances = arrayfun(@(x) calculate_distance(x, equally_spaced_points), equally_spaced_points);

    % Plot results
    figure;

    subplot(3, 1, 1);
    plot(cheb_points, cheb_distances, 'o-', 'LineWidth', 2);
    title('Chebyshev Points');
    xlabel('x_j');
    ylabel('Geometric Mean Distance');

    subplot(3, 1, 2);
    plot(legendre_points, legendre_distances, 'o-', 'LineWidth', 2);
    title('Legendre Points');
    xlabel('x_j');
    ylabel('Geometric Mean Distance');

    subplot(3, 1, 3);
    plot(equally_spaced_points, equally_spaced_distances, 'o-', 'LineWidth', 2);
    title('Equally Spaced Points');
    xlabel('x_j');
    ylabel('Geometric Mean Distance');

    sgtitle(['Geometric Mean Distances for n = ' num2str(n)]);
end

%%Fourier Trig Polynomials function

function P = triginterp(xi,x,y)
% TRIGINTERP Trigonometric interpolation.
% Input:
%   xi  evaluation points for the interpolant (vector)
%   x   equispaced interpolation nodes (vector, length N)
%   y   interpolation values (vector, length N)
% Output:
%   P   values of the trigonometric interpolant (vector)
N = length(x);
% Adjust the spacing of the given independent variable.
h = 2/N;
scale = (x(2)-x(1)) / h;
x = x/scale;  xi = xi/scale;
% Evaluate interpolant.
P = zeros(size(xi));
for k = 1:N
  P = P + y(k)*trigcardinal(xi-x(k),N);
end

function tau = trigcardinal(x,N)
ws = warning('off','MATLAB:divideByZero');
% Form is different for even and odd N.
if rem(N,2)==1   % odd
  tau = sin(N*pi*x/2) ./ (N*sin(pi*x/2));
else             % even
  tau = sin(N*pi*x/2) ./ (N*tan(pi*x/2));
end
warning(ws)
tau(x==0) = 1;     % fix value at x=0

\end{lstlisting}
%%%%%%%%%%%%%%%%%%%%%%%%%%%%%%%%%%%%%%%%%%%%%%%%%%%%%%%%%%%%%%%%%%%%%%%%%%%%

\end{document}